\documentclass[12pt,leqno]{article}
\usepackage{amssymb}
\usepackage{amsmath,amsfonts,amssymb}
\numberwithin{equation}{section} \textheight 240mm \textwidth
160mm \topmargin -5mm \oddsidemargin 0mm
\def\i{\`\i}

\numberwithin{equation}{section}
\newlength{\defbaselineskip}
\DeclareMathAlphabet{\mathpzc}{OT1}{pzc}{m}{n}

\newtheorem{theo}{Theorem}[section]

\newtheorem{lem}{Lemma}[section]

\def\build#1_#2^#3{\mathrel{\mathop{\kern 0pt#1}\limits_{#2}^{#3}}}

\begin{document}

${}$


\vskip 2 cm

\medskip
\centerline{\Large {\bf Equation system describing the radiation
intensity}}

\centerline{\Large {\bf and the air motion with the water phase
transition }}

\medskip
\bigskip
\medskip

\centerline{\bf Meryem Benssaad${}^{1)}$, Hanane Belhireche${}^{2)}$, Steave C.
Selvaduray${}^{3)}$ }

\smallskip

\bigskip
\centerline{${}^{1), 2)}$ University of 8 mai 1945, Guelma,
Alg\'erie}

\centerline{${}^{3)}$ Universit\`a di Torino, Italy}

\medskip

\vskip 0.8 cm

{\bf Abstract.} In this paper we consider the equation system
describing the motion of the air and the variation of the
radiation intensity and the quantity of water droplets in the air,
including also the process of water phase transition. Under a
suitable condition we prove the existence and uniqueness of the
local solution. By eliminating the approximation by
regularization of vapor density and by including the equation of radiation,
this result improves previous ones.  \vskip 0.8 cm

\section{Introduction}

Many attempts of modeling the atmosphere have been made
(most classic ones are for example \cite{[MD]}, \cite{[LTW]}), but it seems that the complete
description of atmospheric phenomena remains to be done.
In \cite{[SF]} the well-posedness of an equation system describing rather
completely the air motion and the water phase transition in the air is shown.
But to obtain this result, the authors have slightly modified the equations,
by replacing the vapor density by its approximation in the equations for the densities
of water vapor, of water droplets and of ice crystals.

In this paper we consider the equation system made up of the equations considered
in  \cite{[SF]} and the equation of the radiation intensity.
We prove the existence and the uniqueness of the local solution to this equation system,
without modifying the function representing the vapor density, but we do not consider the equation
for the density of ice crystals, excluding the phase transition to solid state;
in fact, the behavior of this equation would be similar to that
of the equation for the density of water droplets.

From the technical point of view, for the part of the equations concerning the
velocity and the temperature we follow the scheme of \cite{[SF]}, which is largely based
on the classical techniques of Solonnikov (\cite{[So76]} and others). For the part of
the equations concerning the densities of the air, of water vapor and of water
droplets, we use the techniques developed in \cite{[SF]} and \cite{[AS]} and, as is mentioned above,
improve the proof.
For the equation of radiation intensity, we use partially the techniques developed
in \cite{[BE]}, but we introduce also new techniques for the estimation of the solution in
the norms of $L^\infty$ and of $L^2$.

\section{Equation system}

We will consider our equation system in a bounded domain $ \Omega
\subset \mathbb{R}^3$ with a sufficiently regular boundary $ \partial \Omega$. We
denote by $\varrho(t,x)$ ($t \geq 0$, $ x \in \Omega $) the dry
air density, by $\pi(t,x)$ the water vapor density, by
$\sigma(t,m,x)$ ($ m > 0$ is the mass of a droplet) the
density of liquid water,  by $T(t,x)$ the temperature of the air,
by $v(t,x)= (v_1, v_2, v_3)$ the velocity of the air, by $u(t,m,x)=
(u_1, u_2, u_3) $ the velocity of the droplets, and by
$I_{\lambda}(x,q_1)$ the radiation intensity of wavelength
$\lambda$. Moreover, we denote by $\eta$ and $\zeta$ the
viscosity coefficients, by $\kappa$ the thermic conductivity, by
$c_v$ the specific heat at constant volume, and by $L_{gl}$ the latent
heat relative to the gas-liquid transition.

We assume that the pressure is given by
$$
p=R_0(\frac{\varrho}{\mu_a}+\frac{\pi}{\mu_h})T,
$$
where $R_0$, $\mu_a$ and $\mu_h$ are respectively the universal
constant of gases, the average molar mass of the dry air and the
molar mass of $H_2O$. We assume also that the external force is
given by the gradient of a potential $\Phi$.

In order to describe the motion of the air,
the phase transition of water vapor, the variation of the density of water
vapor, the motion of water droplets, the variation of the radiation intensity and their
interactions, based on the fundamental equations given in \cite{[LL]} and \cite{[Liou]}
and their application to the mentioned phenomena (see \cite{[Met]}, \cite{[SF]}, \cite{[BE]}),
we consider the following system of equations:

\begin{equation}\label{eq-rho}
{\partial \varrho \over \partial t } + \nabla \cdot (\varrho v ) =
0,
\end{equation}
\begin{equation}\label{eq-pi}
{\partial \pi \over \partial t } + \nabla \cdot (\pi  v ) = -
H_{gl}(T,\pi,\sigma(m)),
\end{equation}
\begin{equation}\label{eq-sigma}
\frac{\partial \sigma }{\partial t} + \nabla_{(m,x)} \cdot (\sigma
\widetilde U_{4l} (u,T,\pi   ) ) =
\end{equation}
$$
= [ h_{gl} ( T,\pi ; m ) + B_1 ( \sigma ; m )  - g_1 ( m ) [ \pi -
\overline{ \pi}_{vs}(T) ]^ - ] \sigma +
$$
$$
+ g_0 ( m ) [ N^*  - \widetilde N ( \sigma  ) ]^+ [ {\pi - \overline
{\pi}_{vs}(T)} ]^ + + B_2 ( \sigma ;m),
$$
\begin{equation}\label{ex(2.1)}
(\varrho + \pi ) \big( {\partial v \over \partial t } + ( v \cdot
\nabla ) v \big) = \eta \Delta  v + \big( \zeta + {\eta \over 3 }
\big) \nabla ( \nabla \cdot v ) +
\end{equation}
$$
- R_0 \nabla ( ( {\varrho \over \mu_a} + {\pi \over \mu_h} ) T )-
\big[ \int_0^\infty
 \sigma \left( m \right) dm + \varrho  + \pi
\big] \nabla \Phi ,
$$
\begin{equation}\label{ex(2.4)}
(\varrho +\pi)c_v \Big( {\partial T \over \partial t } +
\sum_{j=1}^3 v_j {\partial T \over \partial x_j } \Big) = \kappa
\Delta T  - R_0 ( \frac{\varrho }{ \mu_a} + \frac{\pi }{ \mu_h} )
T \nabla \cdot  v  +
\end{equation}
$$
+ \eta \sum_{i,j=1}^3 \Big( {\partial v_i \over \partial x_j} +
{\partial v_j \over \partial x_i} - {2 \over 3} \delta_{ij} \nabla
\cdot v \Big) {\partial v_i \over \partial x_j} + \zeta (\nabla
\cdot v)^2 -\nabla\cdot \mathcal E + L_{gl} H_{gl}  ,
$$
\begin{equation}\label{EQ4}
-(q_1\cdot\nabla)I_{\lambda}(t,x,q_1)=
b_{\lambda}(t,x)I_{\lambda}(t,x,q_1)
-J_{\lambda}(t,x,q_1,I_{\lambda},T),
\end{equation}
where
\begin{equation}\label{def-U4}
\nabla_{(m,x)}=(\partial_m,\partial_{x_1},
\partial_{x_2},\partial_{x_3}), \qquad\widetilde U_{4l} (u,T,\pi )
= ( m h_{gl}, u_1 , u_2 , u_3 ),
\end{equation}
\begin{equation}\label{def-E}
\mathcal E = (\mathcal E_1 ,\mathcal E_2 ,\mathcal
E_3 ), \qquad\mathcal E_{j}(t,x)=\int_{0}^{\infty}\int_{S^2}I_{\lambda}(t,x,q_1 )
q_{1j}dq_1d\lambda, \qquad j = 1,2,3,
\end{equation}
and $H_{gl} $, $ h_{gl} $,  $B_1 ( \sigma ; m )$, $B_2 ( \sigma
; m )$, $J_{\lambda} (t,x,q_1, I_{\lambda}, T)$, $
b_{\lambda}(t,x) $, $g_1(m)$, $g_0(m)$, $N^*$, $\widetilde{N}(\sigma)$
are the functions (or numbers) which we are going to precise below. We consider
this equation system for $t\geq0$, $x\in\Omega \subset
\mathbb{R}^3$, $m>0$ and $q_1 \in S^2 = \left\{ {q \in \mathbb{R}^3:|q| = 1} \right\}
$.
The function $I_{\lambda}(t,x,q_1)$ appearing in \eqref{EQ4} and \eqref{def-E}
depends on $t$, but the r\^{o}le of $t$ is only that of parameter. So in the
sequel we write simply $I_{\lambda}(x,q_1)$.

Concerning $H_{gl} $ and $ h_{gl} $, which represent the quantity of condensation
(or evaporation) on all droplets and that on droplets of mass $m$, they can
have some general form. But to fix the idea, we consider $H_{gl}$ and $h_{gl}$
having the form proposed in the modeling \cite{[SF]}, that is

\begin{equation}\label{H-gl}
H_{gl} ( T,\pi , \sigma(\cdot) ) = K_1 \int_0^{\infty}\frac{{S_l (
m )}}{m}\sigma(m)dm \big( {\pi - \overline{ \pi}_{vs}(T) } \big),
\end{equation}
\begin{equation}\label{def-h_gl}
h_{gl}=h_{gl}(T,\pi ,m)=K_{1}\frac{S_{l}(m)}{m}(\pi -\overline{\pi
}_{vs}(T)),
\end{equation}
where $K_1$ is a positive constant, $\overline{ \pi}_{vs}(T)$ is
the density of saturated vapor with respect to the liquid state
and $S_l ( m )$ represents the surface area of the droplet of mass
$m$; for $S_l ( m)$ we suppose that
\begin{equation}\label{eq-S-1}
S_l ( \cdot)  \in C^2 (\mathbf{R}_+), \qquad S_l ( \cdot) \geq 0 ,
\end{equation}
\begin{equation}\label{eq-S-2}
S_l ( m)  = 0 \quad \mbox{for} \ \ 0 \leq m \leq \frac{\overline
{m}_a}{2} \quad (0<\overline{m}_a) ,
\end{equation}
\begin{equation}\label{eq-S-3}
S_l ( m) = c_l m^{2/3} , \ \ \quad \mbox{for } \ \  m \geq \overline
{m}_A \quad (\overline {m}_a<\overline {m}_A<\infty);
\end{equation}
$\overline {m}_a $ and $\overline {m}_A $ should represent the lower and
the upper bounds of the aerosols mass.

The terms $ B_1(\sigma;m)$ and $B_2(\sigma;m)$ are defined by
$$
B_1(\sigma;m)=-m \sigma ( m )\int_0^{\infty } {\beta ( {m,m'}
)\sigma ( {m'})dm'} ,
$$
$$
B_2(\sigma;m)=\frac{m}{2}\int_0^m \beta  ( {m - m',m'} ) \sigma
( {m'} )\sigma ( {m - m'} )dm' ,
$$
$$
\beta  ( {m_1,m_2} )=\beta  ( {m_2,m_1} )
\geq 0 \quad \forall m_1,m_2\in\mathbb{R}_+
$$
where $\beta(m,m')$ denotes the rate of occurrence of coagulation of two droplets with mass
$m$ and  with mass $m'$.
Here we
suppose that for some $M>\overline{ m}_a$,
\begin{equation}\label{(2.25)}
\beta (m', m'')  = 0 \quad \mbox{for} \ \ m'+m'' \geq M.
\end{equation}

The appearance of droplets of mass $m$ is represented by
$ g_0 ( m ) [ N^*  - \widetilde N ( \sigma  ) ]^+ [ {\pi - \overline
{\pi}_{vs(l)}(T)} ]^ + $, where
$N^*$ and $\tilde N ( \sigma  )$ should represent respectively
the total number of aerosols susceptible to the formation of droplets
in the unit volume and the number in the unit volume of
aerosols already present in droplets. The disappearance of droplet
of mass $m$ is represented by $g_1 ( m ) [ \pi -
\overline{ \pi}_{vs(l)}(T) ]^ - \sigma$. For the coefficient
functions  $g_0(m)$ and $g_1(m)$, we suppose that they are
sufficiently regular and
\begin{equation}\label{g0g1g2}
{\rm supp} \, g_0 ( \cdot ) \subset [\overline{m}_a , \overline m_A], \qquad
{\rm supp} \, g_1 ( \cdot ) \subset [0 , \overline{ m}_A].
\end{equation}

The functions $b_{\lambda}(t,x)$ and $J_{\lambda} (t,x,q_1,
I_{\lambda},T)$ are defined by the relations
\begin{equation}\label{def-b-lambda}
b_{\lambda}(t,x)=(a_{\lambda}^{(1)}+r_{\lambda}^{(1)})\varrho(t,x)+
(a_{\lambda}^{(2)}+r_{\lambda}^{(2)})\pi(t,x)
+
\end{equation}
$$
+\int^{\infty}_0(a_{\lambda}^{(3)}(m)+r_{\lambda}^{(3)}(m))\sigma(t,m,x)dm,
$$
\begin{equation}\label{def-J-lambda}
J_{\lambda}(t,x,q_1,I_{\lambda},T)= \frac{1}{4\pi}
r_{\lambda}^{(1)}\varrho(t,x)
\int_{S^2}I_{\lambda}(x,q'_1)P^{(1)}_{\lambda}(q'_1\cdot q_1)dq'_1+
\end{equation}
$$
+\frac{1}{4\pi}r_{\lambda}^{(2)}\pi(t,x)
\int_{S^2}I_{\lambda}(x,q'_1)P^{(2)}_{\lambda}(q'_1\cdot q_1)dq'_1+
$$
$$
+\frac{1}{4\pi}\int_0^{\infty}r^{(3)}_{\lambda}(m)\sigma(t,m,x)dm\int_{S^2}
I_{\lambda}(x,q'_1)P^{(3)}_{\lambda}(q'_1\cdot q_1)dq'_1+
$$
$$
+\big(a_{\lambda}^{(1)}\varrho(t,x)+a_{\lambda}^{(2)}\pi(t,x)+
\int^{\infty}_0a_{\lambda}^{(3)}(m)\sigma(t,m,x)dm\big)B[\lambda,T(t,x)],
$$
where $a_{\lambda}^{(1)}$ is the absorbtion coefficient,
$r_{\lambda}^{(1)}$ is the diffusion coefficient of radiation and
$P_{\lambda}^{(1)}(q'_1\cdot q_1)$ is the diffusion coefficient of
radiation by the dry air of direction $q'_1$ in the direction
$q_1$, while $a_{\lambda}^{(2)}$, $r_{\lambda}^{(2)}$ and
$p_{\lambda}^{(2)}(q'_1\cdot q_1)$ are respectively the absorption
coefficient and the diffusion coefficient of radiation and the
diffusion coefficient of radiation from the direction $q'_1$ to the
direction $q_1$ by the water vapor, and $a_{\lambda}^{(3)}(m)$,
$r_{\lambda}^{(3)}(m)$ and $p_{\lambda}^{(3)}(q'_1\cdot q_1)$ are
respectively the absorption coefficient and the diffusion
coefficient of radiation and the diffusion coefficient of radiation
from the direction $q'_1$ to the direction $q_1$ by liquid water
(or solid water);
$\big(a_{\lambda}^{(1)}\varrho(t,x)+a_{\lambda}^{(2)}\pi(t,x)+
\int^{\infty}_0
a_{\lambda}^{(3)}(m)\sigma(m,t,x)dm\big)B[\lambda,T(x)] $ is the
emission of radiation, where the function $B[\lambda,T]$, called
{\it Planck's function}, is given by
\begin{equation}\label{Plank}
B[\lambda, T ]=\frac{2\pi c^2h}{\lambda^5}(e^{\frac{ch}{k\lambda
T}}-1)^{-1}
\end{equation}
(here $c$ is the speed of light, $h$ is the Planck constant and $k$
is the Boltzmann constant).

For the velocity of  water droplets,  $u(m,t,x)$,
we assume the relation
\begin{equation}
u ( t,m,x) = v ( t,x ) - \frac{1}{{\alpha _l ( m )}}\nabla \Phi ,
\end{equation}
where ${\alpha_l ( m )}$ is the friction coefficient between the
droplets with mass $m$ and the air.

\section{Position of the problem}

Consider the equation system
\eqref{eq-rho}-\eqref{EQ4} in a bounded domain $\Omega \subset
\mathbb{R}^3$ with the initial and boundary conditions
\begin{equation}\label{cond-varrho}
\varrho(0,\cdot)=\varrho_0(\cdot)\in W_p^1(\Omega), \quad
\inf_{x\in\Omega}\varrho_0(x)>0,
\end{equation}
\begin{equation}\label{cond-pi}
\pi(0,\cdot)=\pi_0(\cdot)\in W_p^1(\Omega), \quad
\inf_{x\in\Omega}\pi_0(x)>0,
\end{equation}
\begin{equation}\label{cond-sigma}
\sigma(0,\cdot,\cdot)=\sigma_0(\cdot,\cdot)\in
W_p^1(\mathbb{R}_+\times\Omega),\quad \sigma_0(\cdot,\cdot)\geq 0,
\end{equation}
\begin{equation}\label{cond-v}
v|_{\partial\Omega}=0,\qquad v(0,\cdot)=v_0(\cdot)\in
W_p^{2-\frac{2}{p}},\quad v_0|_{\partial\Omega}=0,
\end{equation}
\begin{equation}\label{cond-T}
\nabla T \cdot n\big|_{\partial\Omega} = 0,\qquad
T(0,\cdot)=T_0(\cdot)\in W_q^{2-\frac{2}{q}}.
\end{equation}

For $\sigma_0$ we suppose

\begin{equation}\label{4.15-sigma-nu-bis}
\exists \overline{ M}' \geq \overline{ m}_A>0 \ \ \mbox{such that} \ \ \sigma_0
(m,\cdot) =
 0 \quad \mbox{if} \ \ m  \in \,
]0, \overline{ m}_a ] \cup [ \overline{ M}' , \infty[,
\end{equation}
\begin{equation}\label{cond-sigma-1}
 \partial_m\sigma_0\in
L^{\infty}(\mathbb{R}_+\times\Omega),
\end{equation}
we define
$$
\Omega_M=\,]0,M[\times \Omega, \qquad \forall M>0.
$$

To specify the boundary conditions  for $
\{I_\lambda \}_{ \lambda >0 }$, it is convenient to transform the
equation \eqref{EQ4} into an integral equation, so that we can rewrite
the equation \eqref{EQ4} in the form
\begin{equation}\label{eqqq1}
\frac{d}{d\alpha}I_\lambda (x+\alpha q_1,q_1) = -b_\lambda(t,x)
I_\lambda (x+\alpha q_1,q_1 )+ J_\lambda (t,x+\alpha
q_1,q_1,\varrho,\pi,\sigma,I_\lambda,T).
\end{equation}
For $ (x,q_1) \in \Omega \times S^2 $ we define
\begin{equation}\label{def-alpha0}
\alpha^0_{(x,q_1)} = \inf \{ \, \alpha \in {\mathbb{R}} \, | \,
x+\alpha' q_1 \in \Omega \ \forall\alpha' \in \, ]\alpha ,0 [ \;
\} .
\end{equation}
The equation \eqref{eqqq1} with the condition
\begin{equation}\label{cond-init-0}
I_\lambda (x+ \alpha^0_{(x,q_1)}q_1, q_1 ) = I_\lambda^0 (x+
\alpha^0_{(x,q_1)}q_1, q_1 )
\end{equation}
can be transformed intro the integral equation
\begin{equation}\label{eq1-I}
I_\lambda (x,q_1)=I^0_\lambda (x+ \alpha^0_{(x,q_1)}q_1, q_1 )
e^{-I_b(x,\alpha^0_{(x,q_1)},q_1)} +
\end{equation}
$$
+\frac{r^{(1)}_\lambda}{4\pi}\int_{
\alpha^0_{(x,q_1)}}^{0}\varrho(t,x+\alpha'q_1)\int_{S^2}
P^{(1)}_\lambda (q'_1\cdot q_1)I_\lambda (x+\alpha'
q_1,q'_1)e^{-I_b(x,\alpha',q_1)} dq_1'd\alpha' +
$$
$$
+\frac{r^{(2)}_\lambda}{4\pi}\int_{
\alpha^0_{(x,q_1)}}^{0}\pi(t,x+\alpha'q_1)\int_{S^2}
P^{(2)}_\lambda (q'_1\cdot q_1)I_\lambda (x+\alpha'
q_1,q'_1)e^{-I_b(x,\alpha',q_1)} dq'_1d\alpha' +
$$
$$
+\frac{1}{4\pi}\int_0^{\infty}r^{(3)}_{\lambda}(m)
\int_{\alpha^0_{(x,q_1)}}\sigma(t,m,x+\alpha'q_1)\int_{S^2}
I_{\lambda}(x+\alpha'q_1,q'_1)P^{(3)}_{\lambda}(q'_1\cdot q_1)\times
$$
$$
\times e^{-I_b(x,\alpha',q_1)}dq_1'd\alpha'dm+
$$
$$
+\int_{ \alpha^0_{(x,q_1)}}^{0}\big(a_{\lambda}^{(1)}\varrho
(t,x+\alpha'q_1)+a_{\lambda}^{(2)}\pi(t,x+\alpha'q_1)+
\int^{\infty}_0a_{\lambda}^{(3)}(m)\sigma(t,m,x+\alpha'q_1)dm\big)\times
$$
$$
\times B[\lambda,T(t,x+\alpha'q_1)]e^{-I_b(x,\alpha',q_1)} d\alpha',
$$
where
\begin{equation}\label{def-I-b}
I_b(x,\alpha,q_1)=\int^0_{\alpha}b_{\lambda}(t,x+\alpha'q_1)d\alpha'.
\end{equation}

We remark that in \eqref{eq1-I} $(x+\alpha^0_{(x,q_1)}q_1,q_1)$
should belong to the set
\begin{equation}\label{def-Xi}
\Xi = \bigcup_{x^0 \in \partial \Omega} \big( \{ x^0 \} \times
{S_-^2}(x^0) \big) ,
\end{equation}
where
\begin{equation}\label{def-S2-minus}
{S_-^2}(x^0)=\{q_1\in S^2\, | \, \exists\, \varepsilon>0\, , \,
x^0+\alpha q_1 \in \Omega, \; \forall \alpha \in\, ]0,\varepsilon[
\, \}
\end{equation}
($x^0 \in \partial\Omega $).

For the diffusion rate $ P^{(i)}_\lambda (q'_1\cdot q_1) $ and for
the diffusion coefficients $r^{(i)}_{\lambda}$  ($i = 1, 2,3 $) we
suppose
\begin{equation}\label{cond0-P}
P^{(i)}_\lambda (q'_1\cdot q_1)\geq 0 \quad \forall q'_1, q_1 \in S^2 ,\qquad
 \frac{1}{4\pi}\int_{S^2}P^{(i)}_\lambda (q'_1\cdot q_1)dq_1=1 \quad
\forall q'_1 \in S^2.
\end{equation}
\begin{equation}\label{hyp-r3-sigma-1}
\sup_{x \in \Omega } \int_0^{\infty} r_{\lambda}^{(3)}(m)
\sigma_0(x , m)dm \leq 4 ,\qquad \int_0^{\infty}
(a_{\lambda}^{(3)}(m)+r_{\lambda}^{(3)}(m))dm\big)\leq \infty,
\end{equation}
moreover, we assume that there is a strictly positive constant
$\varepsilon_1$ such that
\begin{equation}\label{hyp-r1-rho-r2-pi-1} \sup_{x \in \Omega , -1
\leq c \leq 1} \big(r_{\lambda}^{(1)} P^{(1)}_{\lambda}(c)
\varrho_0(x )+ r_{\lambda}^{(2)} P^{(2)}_{\lambda}(c) \pi_0(x )
\big) \leq \frac{\varepsilon_1}{2} ,
\end{equation}
\begin{equation}\label{cond-Kb-P3-epsln1-1}
\sup_{\lambda\in\mathbb{R}_+}( K_b { \sup_{-1 \leq c \leq 1} P^{(3)}_{\lambda} (c ) })^{1/2} +
\frac{{\varepsilon }_1 }{2}< 1 ,
\end{equation}
where
\begin{equation}\label{def-Kb}
K_b = \sup_{ x \in \Omega , q_1 \in S^2 } (  1 -  e^{- 2
I_{b^0}(x, \alpha^0_{(x,q_1)} ,q_1)} ),
\end{equation}
with $I_{b^0}(x,\alpha^0_{(x,q_1)},q_1)$ defined in an analogous way
to \eqref{def-I-b} but with
$$
b^0_{\lambda}(x)=2(a_{\lambda}^{(1)}+r_{\lambda}^{(1)})\varrho_0(x)+
2(a_{\lambda}^{(2)}+r_{\lambda}^{(2)})\pi_0(x)+
2\int^{\infty}_0(a_{\lambda}^{(3)}+r_{\lambda}^{(3)})\sigma_0(m,x)dm+\varepsilon_2
$$
instead of $b_{\lambda}$,

where $ \varepsilon_2 $ is a strictly positive constant and
sufficiently small.

\medskip

It is not restrictive to suppose that the diameter of the domain
$\Omega$ is equal to 1 because we can transform a generic bounded
domain into a domain with diameter equal to 1 by a simple change
of variables.

For the function $\Phi$ we suppose
\begin{equation}\label{4.15-geo-a}
\Phi \in C^3 (\Omega), \qquad \nabla \Phi  \cdot n = 0 \quad
{\hbox{on}} \ \ \partial \Omega
\end{equation}
({$n$ is the unit outward normal vector to $\partial \Omega$}).

\bigskip

The main result of this paper is the following
theorem.

\begin{theo}\label{maintheorem}

Let us assume $p> 4$, $2q > p> q >3 $ and the conditions
\eqref{cond-varrho}-\eqref{cond-T}, \eqref{cond-init-0},
\eqref{cond0-P}-\eqref{cond-Kb-P3-epsln1-1}.
Then there exists a number $\overline t >0$ such that the problem
\eqref{eq-rho}-\eqref{EQ4} admits a unique solution $(\varrho, \pi,
\sigma, v, T, I_{\lambda})$ with the following properties :
\begin{equation}\label{ex(4.21)}
\varrho \in C^0 ([0, \overline t ];W^1_p(\Omega)) , \qquad
\inf_{(t,x) \in Q_{\overline t} }\varrho (t,x) >0 ,
\end{equation}
\begin{equation}\label{ex(4.22)}
\pi \in C^0([0,\overline t ]; W^1_p(\Omega)) , \quad
\pi \geq 0 ,
\end{equation}
\begin{equation}\label{ex(4.23)}
\sigma \in C^0 ([0, \overline t ];
W^1_p( \mathbf{R}_+ \times \Omega ) ) ,
\quad \sigma \geq 0 ,\quad {\partial _m \sigma
\in C^0 ([0, \overline t ] ;L^\infty  (\mathbf{R}_+ \times \Omega ))},
\end{equation}

\begin{equation}\label{ex(4.20)}
v \in W^{2,1}_p([0, \overline t ] \times \Omega   ) , \qquad
T \in W^{2,1}_q([0, \overline t ] \times \Omega ) , \quad
T > 0 ,
\end{equation}

\begin{equation}
I_{\lambda}\in L^{\infty}(\Omega\times S^2).
\end{equation}

\end{theo}

\section{Equation of radiation intensity}

In this section, supposing that $\varrho$, $\pi$, $\sigma$ and $T$
are given, we prove the existence and
uniqueness of the solution to the equation \eqref{eq1-I} with
fixed $\lambda$ and $t$, and give an estimate of the difference of two solutions of
\eqref{eq1-I} with the same  $\lambda$ and $t$ but different
$\varrho$, $\pi$, $\sigma$ and $T$.

\begin{lem}\label{lemme01}
Let be $I_\lambda^0 (x^0 , q_1 ) $ a non-negative measurable function
defined on $\Xi $. We suppose that
\begin{equation}\label{cond-I-0}
\sup_{ (x^0 , q_1) \in \Xi } I_\lambda^0 (x^0,q_1) < \infty .
\end{equation}
If the functions $\varrho(t, \cdot )$, $\pi(t, \cdot )$ and $T(t,
\cdot )$ are given in ${L}^{\infty}(\Omega )$ and $\sigma(t,\cdot
, \cdot)$ is given in ${L}^{\infty}(\mathbb{R}_+\times\Omega )$,
then the equation \eqref{eq1-I} admits a unique solution
$I_\lambda $ in ${L}^{\infty}(\Omega\times S^2)$ and we have
\begin{equation}\label{sup-I-lambda}
\sup_{(x,q_1)\in\Omega\times S^2}I_{\lambda}(x,q_1)\leq
\frac{1}{\varepsilon_b}\big[\sup_{(x^0,q_1)\in\Xi}I_{\lambda}^{(0)}(x^0,q_1)+
\sup_{\frac{\overline T_0^{(-)}}{2}\leq T
\leq\frac{3\overline T_0^{(+)}}{2}}B[\lambda,T] \big],
\end{equation}
where
$$
\varepsilon_b=\varepsilon_b(t)=\inf_{(x,q_1)\in\Omega\times S^2}
e^{-I_b(x,\alpha^0_{x,q_1},q_1)},$$
$$ \overline T_0^{(-)}=\inf_{x\in\Omega}T_0(x),\qquad
\overline T_0^{(+)}=\sup_{x\in\Omega}T_0(x).
$$

\end{lem}

{\bf Proof.} First, supposing that $I_{\lambda}(x,q_1)$
is a solution to \eqref{eq1-I}, we prove the inequality
\eqref{sup-I-lambda}, setting
$$
\overline A=\sup_{(x,q_1)\in\Omega\times S^2}I_{\lambda}(x,q_1),
\qquad \overline B=\sup_{\frac{\overline T_0^{(-)}}{2}\leq T
\leq\frac{3\overline T_0^{(+)}}{2}}B[\lambda,T] ,
$$
$$
\overline I=\sup_{(x^0,q_1)\in\Xi}I_{\lambda}^{(0)}(x^0,q_1).
$$
From the equation \eqref{eq1-I}, we deduce
$$
I_{\lambda}(x,q_1)\leq \overline I+\overline A\int^{0}_{\alpha^0_{(x,q_1)}}
b_{\lambda}(t,x+\alpha'q_1)e^{-I_b(x,\alpha',q_1)}d\alpha'+
$$
$$
+\overline B\int^{0}_{\alpha^0_{(x,q_1)}}
b_{\lambda}(t,x+\alpha'q_1)e^{-I_b(x,\alpha',q_1)}d\alpha'.
$$

As $$\int^{0}_{\alpha^0_{(x,q_1)}}
b_{\lambda}(t,x+\alpha'q_1)e^{-I_b(x,\alpha',q_1)}d\alpha'<1-\varepsilon_b$$
(see \eqref{def-b-lambda}), we obtain

$$
\overline A\leq\overline I+(1-\varepsilon_b)(\overline A+\overline B).
$$
From this inequality we deduce \eqref{sup-I-lambda}.

To prove the existence and the uniqueness  of the solution $I_{\lambda}$, we
denote by $G(I_{\lambda})$ the second member of \eqref{eq1-I}. So
we have
\begin{equation}
|G(I_{\lambda}^{[1]})(x,q_1)-G(I_{\lambda}^{[2]})(x,q_1)|\leq
\end{equation}
$$
\leq\sup_{(x,q_1)\in \Omega\times S^2}|I_{\lambda}^{[1]}-I_{\lambda}^{[2]}|\frac{1}{4\pi}
\int_{\alpha^0_{(x,q_1)}}\int_{S^2}\big(r^{(1)}_{\lambda}\varrho(t,x+\alpha'q_1)
P^{(1)}_\lambda (q'_1,q_1)+
$$
$$
+r^{(2)}_{\lambda}\pi(t,x+\alpha'q_1)P^{(2)}_\lambda (q'_1,q_1)
+\int_0^{\infty}r^{(3)}_{\lambda}(m)
\sigma(t,m,x+\alpha'q_1)P^{(3)}_{\lambda}(m,q'_1,q_1)dm\big)\times
$$
$$
\times e^{-I_b(x,\alpha',q_1)}dq_1'd\alpha'.
$$
Hence, using the properties on $ P^{(i)}_\lambda $ (see
\eqref{cond0-P}), we obtain
\begin{equation}\label{oper-G}
|G(I_{\lambda}^{[1]})(x,q_1)-G(I_{\lambda}^{[2]})(x,q_1)|\leq
\end{equation}
$$
\leq\int^0_{\alpha^0_{(x,q_1)}}b_{\lambda}(t,x+\alpha'q_1)
e^{-I_b(x,\alpha',q_1)}d\alpha'
\sup_{(x,q_1)\in \Omega\times S^2}|I_{\lambda}^{[1]}-I_{\lambda}^{[2]}|\leq
$$
$$
\leq\big(1-e^{-I_b(x,\alpha^0_{x,q_1},q_1)}\big)
\sup_{(x,q_1)\in \Omega\times S^2}|I_{\lambda}^{[1]}-I_{\lambda}^{[2]}|,
$$
where $I_b(x,\alpha',q_1)$ is defined by
\eqref{def-I-b}. That is, the operator $G(\cdot)$ is a contraction in
$L^{\infty}(\Omega\times S^2)$, so that the equation \eqref{eq1-I}
has a unique solution in $L^{\infty}(\Omega\times S^2)$. $\quad
\square$

\medskip
\medskip

Now, we are going to prove an estimate for the divergence of $\mathcal E$.

\begin{lem}\label{lemme03}
Let us assume $\varrho(t,\cdot),\pi(t,\cdot)\in L^{p}(\Omega)$,
$\sigma(t,\cdot,\cdot)\in L^p(L^{\infty}(\mathbb{R}_+);\Omega)$,
$I_{\lambda}\in L^{\infty}( \Omega\times S^2)$,
$B[\lambda,T(\cdot)]\in L^{\infty}(\Omega)$.
Then there exists a positive constant $c$ such that

\begin{equation}\label{estim-E-p}
\|\nabla\cdot \mathcal E\|_{L^q}\leq
c\big(\big\|\int_0^{\infty}I_{\lambda}(x,q)d\lambda\big\|_{L^{\infty}
(\Omega\times S^2)}+
\big\|\int_0^{\infty}B[\lambda,T(x)]d\lambda\big\|_{L^{\infty}(\Omega)}\big)\times
\end{equation}
$$
\times\big(\|\varrho\|_{L^p(\Omega)}+\|\pi\|_{L^p(\Omega)}+
\|\sigma\|_{L^p(\Omega;L^{\infty}(\mathbb{R}_+))}\big).
$$
\end{lem}
{\bf Proof.} It is not difficult to obtain, by elementary calculus (see also the conditions
\eqref{cond0-P}, \eqref{hyp-r3-sigma-1}), the inequality \eqref{estim-E-p} from the
definition of $\mathcal E(\cdot)$ (see \eqref{EQ4}, \eqref{def-E}). \ $\square$

\medskip

We prove also some estimates for the
difference between two possible functions representing the
radiation intensity $I_\lambda$.

Let $I_\lambda^{[i]} (x,q_1 )$, $i = 1, 2$, two functions in $
L^\infty (\Omega \times S^2)$ verifying the equation \eqref{eq1-I}
with $\varrho=\varrho_i$, $\pi=\pi_i$, $\sigma=\sigma_i$, $T=T_i$,
$i=1,2$. Then we have
\begin{equation}\label{eq-diff-I-lam-estimL2}
I_{\lambda}^{[1] }(x,q_1 ) -I_{\lambda}^{[2] }(x,q_1 ) =
 \Delta_0 I_\lambda+\Delta_1 I_\lambda+\Delta_2 I_\lambda ,
\end{equation}
where
\begin{equation}\label{def-Delta0-Ilam}
\Delta_0 I_\lambda=I_{\lambda}^{(0)} (x+ \alpha^{0}_{(x,q_1)}q_1,q_1)
(e^{-I_{b_1}(x,\alpha^0,q_1)} -e^{-I_{b_2}(x,\alpha^0,q_1)}) +
\end{equation}
$$
+\frac{1}{4\pi}\int_{\alpha^0_{(x,q_1)}}^0\int_{S^2} \big(
r_{\lambda}^{(1)} P^{(1)}_{\lambda}(q_1' \cdot q_1) (\varrho_1-
\varrho_2)(x+\alpha'q_1)+ r_{\lambda}^{(2)} P^{(2)}_{\lambda}(q_1'
\cdot q_1) (\pi_1-\pi_2)(x+\alpha'q_1)+
$$
$$
+\int_0^{\infty}r_{\lambda}^{(3)}(m)P^{(3)}_{\lambda}(q_1' \cdot
q_1 ) (\sigma_1-\sigma_2)(m,x+\alpha'q_1)dm \big)\times
$$
$$
\times e^{-I_{b_2}(x,\alpha',q_1)}I_{\lambda}^{[2] }
(x+\alpha'q_1, q'_1)dq'_1d\alpha'+
$$
$$
+\frac{1}{4\pi}\int_{\alpha^0_{(x,q_1)}}^0\int_{S^2}
\big(r_{\lambda}^{(1)}P^{(1)}_{\lambda}(q_1' \cdot
q_1)\varrho_1(x+\alpha'q_1)+ r_{\lambda}^{(2)}
P^{(2)}_{\lambda}(q_1' \cdot q_1)\pi_1(x+\alpha'q_1)+
$$
$$
+\int_0^{\infty}r_{\lambda}^{(3)}(m)P^{(3)}_{\lambda}(q_1' \cdot
q_1 )\sigma_1(m,x+\alpha'q_1)dm \big)\times
$$
$$
\times(e^{-I_{b_1}(x,\alpha',q_1)}-e^{-I_{b_2}(x,\alpha',q_1)})
I_{\lambda}^{[2] } (x+\alpha'q_1,q'_1)dq'_1d\alpha'+
$$
$$
+\int_{\alpha^0_{(x,q_1)}}^0\big(a_{\lambda}^{(1)} \varrho_1
(x+\alpha'q_1) +a_{\lambda}^{(2)}\pi_1(x+\alpha'q_1)+
\int_0^{\infty} a_{\lambda}^{(3)}(m) \sigma_1 (m,x+\alpha'q_1)
\big)\times
$$
$$
\times\big[(B[\lambda,T_1]-B[\lambda,T_2])+
(e^{-I_{b_1}(x,\alpha',q_1)}-e^{-I_{b_2}(x,\alpha',q_1)})B[\lambda,T_2]\big]d\alpha'+
$$
$$
+\int_{\alpha^0_{(x,q_1)}}^0 \big( a_{\lambda}^{(1)}
(\varrho_1-\varrho_2) (x+\alpha'q_1) + a_{\lambda}^{(2)}
(\pi_1-\pi_2)(x+\alpha'q_1)+
$$
$$
+ \int_0^{\infty} a_{\lambda}^{(3)}(m) (\sigma_1 - \sigma_2)
(m,x+\alpha'q_1)dm\big) e^{-I_{b_2}(x,\alpha',q_1)}
B[\lambda,T_2]d\alpha' ,
$$

\begin{equation}\label{def-Delta1-Ilam}
\Delta_1 I_\lambda = \frac{1}{4\pi} \int_{\alpha^0_{(x,q_1)}}^0
\int_{S^2} \big(r_{\lambda}^{(1)}P^{(1)}_{\lambda}(q_1' \cdot q_1)
\varrho_1(x +\alpha'q_1)+
\end{equation}
$$
+r_{\lambda}^{(2)} P^{(2)}_{\lambda}(q_1'
\cdot q_1) \pi_1(x+\alpha'q_1) \big) \times
$$
$$
\times e^{-I_{b_1}(x,\alpha',q_1)} (I_{\lambda}^{[1]} -
I_{\lambda}^{[2] }) (x+\alpha'q_1, q'_1)dq'_1d\alpha' ,
$$

\begin{equation}\label{def-Delta2-Ilam}
\Delta_2 I_\lambda = \frac{1}{4\pi} \int_{\alpha^0_{(x,q_1)}}^0
\int_{S^2} \int_0^{\infty}r_{\lambda}^{(3)}(m) P^{(3)}_{\lambda}
(q_1' \cdot q_1 ) \sigma_1(m,x+\alpha'q_1)dm \times
\end{equation}
$$
\times e^{-I_{b_1}(x,\alpha',q_1)} (I_{\lambda}^{[1]} -
I_{\lambda}^{[2] }) (x+\alpha'q_1, q'_1)dq'_1d\alpha' ;
$$
here $ I_{b_i}(x,\alpha',q_1)$ denotes the
function $ I_b (x, \alpha,q_1) =\int_{\alpha}^0 b_{\lambda} (t,x+
\alpha'q_1) d\alpha' $ with $ \varrho = \varrho_i$, $ \pi =
\pi_i$, $ \sigma = \sigma_i$, $i = 1, 2$.

\medskip

\begin{lem}\label{lemme-9-1}
Let us assume \eqref{hyp-r3-sigma-1} and \eqref{hyp-r1-rho-r2-pi-1}.
Then we have
\begin{equation}\label{estim-Delta1-Ilam}
| \Delta_1 I_\lambda | \leq \frac{\varepsilon_1}{4\pi}
\int_{\alpha^0_{(x,q_1)}}^0 \int_{S^2} |I_{\lambda}^{[1]} -
I_{\lambda}^{[2] } | (x+\alpha'q_1, q'_1) dq'_1d\alpha' ,
\end{equation}
\begin{equation}\label{estim-Delta2-Ilam}
| \Delta_2 I_\lambda | \leq  K_b^{1/2} \Big( \frac{1}{4\pi}
\int_{S^2} P^{(3)}_{\lambda} (q_1' \cdot q_1 )
\int_{\alpha^0_{(x,q_1)}}^0 (I_{\lambda}^{[1]} - I_{\lambda}^{[2]
})^2 (x+\alpha'q_1, q'_1) d\alpha ' d q'_1 \Big)^{1/2} ,
\end{equation}
where
\begin{equation}\label{def-Kb}
K_b = \sup_{ x \in \Omega , q_1 \in S^2 } (  1 -  e^{- 2
I_{b_1}(x, \alpha^0_{(x,q_1)} ,q_1)} ) .
\end{equation}
\end{lem}

\medskip

{\bf Proof}. \ The inequality \eqref{estim-Delta1-Ilam} results
immediately from \eqref{hyp-r1-rho-r2-pi-1}.

On the other hand, using the Cauchy-Schwartz inequality, we have
$$
| \Delta_2 I_\lambda | \leq
$$
$$
\leq \frac{ 1}{4\pi} \int_{S^2} P^{(3)}_{\lambda} (q_1' \cdot q_1
) \Big( \int_{\alpha^0_{(x,q_1)}}^0 \big( \int_0^{\infty}
r_{\lambda}^{(3)}(m) \sigma_1(m,x+\alpha'q_1)dm \big)^2 e^{- 2
I_{b_1}(x,\alpha',q_1)} d\alpha ' \Big)^{1/2} \times
$$
$$
\times \Big( \int_{\alpha^0_{(x,q_1)}}^0 (I_{\lambda}^{[1]} -
I_{\lambda}^{[2] })^2 (x+\alpha'q_1, q'_1) d\alpha ' \Big)^{1/2} d
q'_1 .
$$
We remark that the condition \eqref{hyp-r3-sigma-1} and the
definitions of $ I_{b_1}(x,\alpha',q_1) $ and  $ b_{\lambda} (t,
x) $ imply
$$
\int_{\alpha^0_{(x,q_1)}}^0 \big( \int_0^{\infty}
r_{\lambda}^{(3)}(m) \sigma_1(m,x+\alpha'q_1)dm \big)^2 e^{- 2
I_{b_1}(x,\alpha',q_1)} d\alpha ' \leq
$$
$$
\leq \int_{\alpha^0_{(x,q_1)}}^0 ( 2 \int_0^{\infty}
r_{\lambda}^{(3)}(m) \sigma_1(m,x+\alpha'q_1)dm ) e^{- 2
I_{b_1}(x,\alpha',q_1)} d\alpha ' \leq
$$
$$
\leq 1 -  e^{- 2 I_{b_1}(x, \alpha^0_{(x,q_1)} ,q_1)} .
$$
So we deduce
$$
| \Delta_2 I_\lambda | \leq \frac{ K_b^{1/2}}{4\pi} \int_{S^2}
P^{(3)}_{\lambda} (q_1' \cdot q_1 ) \Big(
\int_{\alpha^0_{(x,q_1)}}^0 (I_{\lambda}^{[1]} - I_{\lambda}^{[2]
})^2 (x+\alpha'q_1, q'_1) d\alpha ' \Big)^{1/2} d q'_1 .
$$
As $ \frac{ 1}{4\pi} \int_{S^2} P^{(3)}_{\lambda} (q_1' \cdot q_1
) d q'_1 = 1 $, by applying again the Cauchy-Schwartz inequality
to the right-hand side of the last inequality, we obtain
\eqref{estim-Delta2-Ilam}. $ \square$

\medskip

\begin{lem}\label{lemme-9-2}
Let $ \varphi (x) \geq 0$ for all $x\in \Omega$. Then we have
\begin{equation}\label{int-multiple}
\frac{1}{4\pi} \int_\Omega \int_{S^2} \int_{\alpha^0_{(x,q_1)}}^0
\varphi (x + \alpha' q_1 ) d\alpha' dq_1 dx \leq \int_\Omega
\varphi (x ) dx .
\end{equation}

\end{lem}

\medskip

{\bf Proof}. \ From
$$
\int_{S^2} \int_{\alpha^0_{(x,q_1)}}^0 \varphi (x + \alpha' q_1 )
d\alpha' dq_1 \leq \int_\Omega \varphi (x') \frac{1}{ | x - x'|^2
} dx'
$$
we obtain
$$
\frac{1}{4\pi} \int_\Omega \int_{S^2} \int_{\alpha^0_{(x,q_1)}}^0
\varphi (x + \alpha' q_1 ) d\alpha' dq_1 dx \leq \frac{1}{4\pi}
\int_\Omega \Big( \int_\Omega \varphi (x') \frac{1}{ | x - x'|^2 }
dx' \Big) dx =
$$
$$
= \int_\Omega \varphi (x') dx' \Big( \frac{1}{4\pi} \int_\Omega
\frac{1}{ | x - x'|^2 } dx \Big) = \int_\Omega \varphi (x') dx'
\Big( \frac{1}{4\pi} \int_0^{{\rm dist} (x' , \partial \Omega)}
\frac{1}{ r^2 } {\mu}_2 (\Sigma^{x'}_r ) dr \Big) ,
$$
where $ {\mu}_2 (\cdot )$ is the Hausdorf measure of dimension 2
and
$$ \Sigma^{x'}_r = \{ \, x \in \Omega \, | \, |x - x' |=r \, \}
.
$$
As
$$
{\mu}_2 (\Sigma^{x'}_r ) \leq 4\pi r^2 , \qquad {{\rm dist} (x' ,
\partial \Omega)} \leq 1 ,
$$
we deduce \eqref{int-multiple}. \ $ \square $

\medskip

\begin{lem}\label{lemme-9-3}

We assume the condition
\eqref{cond-Kb-P3-epsln1-1}.
Then there exists a constant $C$ such that
\begin{equation}\label{estimL2-diff-Ilam}
\int_0^{\infty}\| I_{\lambda}^{[1] } -I_{\lambda}^{[2] } \|_{L^2 (\Omega \times
S^2 )}^2d\lambda \leq C \big[ \|\varrho_2-\varrho_1 \|^2_{L^2 ( \Omega )} + \|
\pi_2-\pi_1  \|_{L^2 ( \Omega )}^2  +
\end{equation}
$$+ \| \sigma_2-\sigma_1 \|_{L^2 ( \Omega_{\overline M_1}
)}^2 + \|T_2-T_1 \|_{L^2 ( \Omega )}^2 \big] .
$$
\end{lem}

\medskip

{\bf Proof}.
Using \eqref{estim-Delta1-Ilam} and the Cauchy-Schwartz
inequality, we have
$$
\int_\Omega \int_{S^2} \Delta_1 I_\lambda (I_{\lambda}^{[1] }
-I_{\lambda}^{[2] } ) dq_1 dx \leq \frac{\varepsilon_1}{4\pi}
\Big( \int_\Omega \int_{S^2} | I_{\lambda}^{[1] }
-I_{\lambda}^{[2] } |^2 dq_1 dx \Big)^{1/2} \times
$$
$$
\times \Big( \int_\Omega dx \int_{S^2} dq_1 \Big(
\int_{\alpha^0_{(x,q_1)}}^0 \int_{S^2} |I_{\lambda}^{[1]} -
I_{\lambda}^{[2] } | (x+\alpha'q_1, q'_1) dq'_1d\alpha' \Big)^2
\Big)^{1/2} .
$$
Now, taking into account the condition $
|\alpha^0_{(x,q_1)}|\leq 1$, we have
$$
\frac{ 1}{ \sqrt{4\pi}} \int_{\alpha^0_{(x,q_1)}}^0 \int_{S^2}
|I_{\lambda}^{[1]} - I_{\lambda}^{[2] } | (x+\alpha'q_1, q'_1)
dq'_1d\alpha' \leq
$$
$$
\leq \Big( \int_{\alpha^0_{(x,q_1)}}^0 \int_{S^2}
|I_{\lambda}^{[1]} - I_{\lambda}^{[2] } |^2 (x+\alpha'q_1, q'_1)
dq'_1d\alpha' \Big)^{1/2} .
$$
Changing the order of the integration with respect to $ q_1 $ and
with respect to $ q'_1 $ and applying Lemma \ref{lemme-9-2} to
$ \varphi (x + \alpha' q_1 ) = |I_{\lambda}^{[1]} -
I_{\lambda}^{[2] } | (x+\alpha'q_1, q'_1) $ for each fixed $ q'_1
$, we obtain
$$
\frac{ 1}{ (4\pi)^2} \int_\Omega dx \int_{S^2} dq_1 \Big(
\int_{\alpha^0_{(x,q_1)}}^0 \int_{S^2} |I_{\lambda}^{[1]} -
I_{\lambda}^{[2] } | (x+\alpha'q_1, q'_1) dq'_1d\alpha' \Big)^2
\leq
$$
$$
\leq \frac{ 1}{{4\pi}} \int_{S^2} dq'_1 \int_\Omega dx \int_{S^2}
\int_{\alpha^0_{(x,q_1)}}^0 |I_{\lambda}^{[1]} - I_{\lambda}^{[2]
} |^2 (x+\alpha'q_1, q'_1) d\alpha' dq_1 \leq
$$
$$
\leq \int_{S^2} dq'_1 \int_\Omega dx |I_{\lambda}^{[1]} -
I_{\lambda}^{[2] } |^2 (x , q'_1) = \|I_{\lambda}^{[1]} -
I_{\lambda}^{[2] } \|_{ L^2 (\Omega \times S^2) }^2 .
$$
We deduce
\begin{equation}\label{estim-D1I-L2}
\int_\Omega \int_{S^2} \Delta_1 I_\lambda (I_{\lambda}^{[1] }
-I_{\lambda}^{[2] } ) dq_1 dx \leq {\varepsilon_1}
\|I_{\lambda}^{[1]} - I_{\lambda}^{[2] } \|_{ L^2 (\Omega \times
S^2) }^2 .
\end{equation}

On the other hand, according to \eqref{estim-Delta2-Ilam} we have
$$
\int_\Omega \int_{S^2} \Delta_2 I_\lambda (I_{\lambda}^{[1] }
-I_{\lambda}^{[2] } ) dq_1 dx \leq K_b^{1/2} \|I_{\lambda}^{[1]} -
I_{\lambda}^{[2] } \|_{ L^2 (\Omega \times S^2) } \times
$$
$$
\times \Big( \frac{1}{ 4 \pi} \int_\Omega dx \int_{S^2} dq_1
\int_{S^2} P^{(3)}_{\lambda} (q_1' \cdot q_1 )
\int_{\alpha^0_{(x,q_1)}}^0 (I_{\lambda}^{[1]} - I_{\lambda}^{[2]
})^2 (x+\alpha'q_1, q'_1) d\alpha ' d q'_1 \Big)^{1/2} .
$$
From Lemma \ref{lemme-9-2}, we obtain
$$
\frac{1}{ 4 \pi} \int_\Omega dx \int_{S^2} dq_1 \int_{S^2}
P^{(3)}_{\lambda} (q_1' \cdot q_1 ) \int_{\alpha^0_{(x,q_1)}}^0
(I_{\lambda}^{[1]} - I_{\lambda}^{[2] })^2 (x+\alpha'q_1, q'_1)
d\alpha ' d q'_1 \leq
$$
$$
\leq { \sup_{-1 \leq c \leq 1} P^{(3)}_{\lambda} (c ) }
\int_\Omega \int_{ S^2 } (I_{\lambda}^{[1]} - I_{\lambda}^{[2]
})^2 (x , q'_1)  d q'_1 dx = { \sup_{-1 \leq c \leq 1}
P^{(3)}_{\lambda} (c ) } \|I_{\lambda}^{[1]} - I_{\lambda}^{[2] }
\|_{ L^2 (\Omega \times S^2) }^2 .
$$
Consequently, we have
\begin{equation}\label{estim-D2I-L2}
\int_\Omega \int_{S^2} \Delta_2 I_\lambda (I_{\lambda}^{[1] }
-I_{\lambda}^{[2] } ) dq_1 dx \leq ( K_b { \sup_{-1 \leq c \leq 1}
P^{(3)}_{\lambda} (c ) })^{1/2} \|I_{\lambda}^{[1]} -
I_{\lambda}^{[2] } \|_{ L^2 (\Omega \times S^2) }^2 .
\end{equation}

For the term $\Delta_0 I_\lambda$, using the intermediate value
theorem and Lemma \ref{lemme-9-2}, we find
$$
\int_\Omega \int_{S^2} \Delta_0 I_\lambda (I_{\lambda}^{[1] }
-I_{\lambda}^{[2] } ) dq_1 dx \leq C_{\lambda}\big( \|\varrho_2-\varrho_1 \|_{L^2 ( \Omega )} + \|
 \pi_2-\pi_1  \|_{L^2 ( \Omega )}+ $$
$$ + \| \sigma_2-\sigma_1 \|_{L^2 ( \Omega_{\overline M_1}
)} + \|T_2-T_1 \|_{L^2 ( \Omega )} \big)
\| I_{\lambda}^{[1] } -I_{\lambda}^{[2] } \|_{L^2 (\Omega \times
S^2 )},
$$
where $C_{\lambda}$ is a constant depending on $I_{\lambda}^{(0)}$,
 $B[\lambda,T_2]$ and
$\frac{\partial}{\partial T}B[\lambda,T]$. $\square$

\section{Linear equations for densities}

In this section we study the equations
\eqref{eq-rho}-\eqref{eq-sigma} for  $\varrho$, $\pi$, $\sigma$
with a given $( {v,T} ) = ( {\overline v ,\overline T } )$.
In the following proofs, we will write simply $c$ to denote
a constant if its specific value will not be used in the sequel;
of course, constants $c$ in different inequalities are different in general.

We introduce the functions spaces
\begin{equation}\label{Theta-v1}
\Theta_{t_1}^{(v)}=\{ v\in W_p^{2,1} (Q_{t_1})| \, v\,
\mbox{satisfies}\,\eqref{cond-v}\},
\end{equation}
\begin{equation}\label{Theta-v2}
\Theta_{t_1}^{(T)}=\{T\in W_q^{2,1} (Q_{t_1})| \,
T\,\mbox{satisfies}\,\eqref{cond-T}\},
\end{equation}
where $Q_{t_1}=\,]0,t_1[\times \Omega$.

Let be $(\overline v,\overline
T)\in\Theta_{t_1}^{(v)}\times\Theta_{t_1}^{(T)}$. We consider the
equations \eqref{eq-rho}-\eqref{eq-sigma} for $\varrho$, $\pi$,
$\sigma$ with $( {v,T} ) = ( {\overline v ,\overline T } )$
\begin{equation}\label{eq-varrho-lin}
\partial_t{\varrho}+\nabla\cdot(\varrho\overline{v})=0,
\end{equation}
\begin{equation}\label{eq-pi-lin}
\partial_t{\pi}+\nabla\cdot(\pi\overline{v})=-H_{gl}
(\overline{T},\pi,\sigma) ,
\end{equation}
\begin{equation}\label{eq-sigma-lin}
\frac{\partial \sigma }{\partial t} + \nabla_{(m,x)} \cdot (\sigma
\widetilde U_{4l} (\overline u,\overline T,\pi) ) =
\end{equation}
$$
= \big[ h_{gl} ( {\overline T ,\pi ;m} ) + B_1 ( \sigma ;m )
 - g_1 ( m )
[ \pi  - \overline {\pi}_{vs(l)} ( {\overline T } ) ]^- \big]\sigma +
$$
$$
 +
g_0 ( m ) [ N^*  - \widetilde N ( {\sigma  } ) ]^+ [ \pi - \overline
{\pi}_{vs(l)} ( {\overline T } )]^+ + B_2 ( \sigma ; m ),
$$
where
\begin{equation}\label{u-bar-v-bar}
\overline u  = \overline v  - \frac{1}{{\alpha _l \left( m
\right)}} \nabla \Phi .
\end{equation}

Before the complete analysis of the equation system
\eqref{eq-varrho-lin}-\eqref{eq-sigma-lin},
in this section, we study the equation \eqref{eq-varrho-lin}, which is linear, and the
linearized equations for \eqref{eq-pi-lin} and
\eqref{eq-sigma-lin}.

For the equation \eqref{eq-varrho-lin} we know the following
result (see for example \cite{[BFY]}, \cite{[FC]}).

\medskip

\begin{lem}\label{lem-varrho}
Let be $\overline v \in \Theta^{(v)}_{t_1}$. The equation
\eqref{eq-varrho-lin} with the condition \eqref{cond-varrho}
admits a unique solution
$$\varrho \in C^0 ([0, t_1
];W^1_p(\Omega))
$$
and we have
\begin{equation}\label{ex(6.7)}
{\|\varrho(t,\cdot)\|}_{W_{p}^1(\Omega)}^p \leq q_\varrho(t),\qquad 0
< \alpha_{\varrho(t)} \leq \varrho(t,x) \leq \beta_{\varrho(t)} <
\infty \quad \mbox{in} \quad Q_{t_{1}},
\end{equation}
where
\begin{equation}\label{ex(6.8)}
\alpha_\varrho(t)=\inf_{x\in\Omega} \varrho_0(x)
\exp(-cR_{(\overline{v},t)}t^{\frac{p-1}{p}}),
\end{equation}
$$
\beta_\varrho(t)=\sup_{x\in\Omega}\varrho_0(x)
\exp(cR_{(\overline{v},t)}t^{\frac{p-1}{p}}),$$
$$
q_{\varrho(t)}={\|\varrho_0\|}_{W_{p}^1(\Omega)}^p
 \exp(cR_{(\overline{v},t)}t^{\frac{p-1}{p}}),\quad
R_{(\overline{v},t)}={\|\overline{v}\|}_{{W_{p}^{2,1}}(Q_t)}.
$$
\end{lem}

For the linearized equation of \eqref{eq-pi-lin}, we have the
following lemma.

\begin{lem}\label{lem-pi}
Let be $\overline v \in \Theta^{(v)}_{t_1} $, $\overline T \in
\Theta^{(T)}_{t_1}$, $ \overline \pi \in C^0([0,t_1];
W^1_p(\Omega)) $, $ \overline \sigma \in C^0([0,t_1];
W^1_p(\mathbf{R}_+ \times \Omega))$. We suppose that there exists
a positive number $\overline M_1$ such that
$supp(\overline{\sigma}(t,\cdot,\cdot))\subset \Omega_{\overline M_1}$ for every $t \in
[0,t_1]$.
Then the equation
\begin{equation}\label{ex(5-eq-pi-bar2)}
\partial_t{\pi}+\nabla\cdot(\pi.\overline{v})=-H_{gl}(\overline{T},{\overline{\pi}},\overline{\sigma})
\end{equation}
with the initial condition \eqref{cond-pi} admits a unique
solution $ \pi\in\mathcal{C}^0([0,t_1];W_{p}^1(\Omega))$ and we
have
\begin{equation}\label{ex(5.14)}
{\|\pi(t,\cdot)\|}_{W_{p}^1(\Omega)}^p \leq q_\pi(t),
\end{equation}
where
\begin{equation}\label{ex(5.15)}
q_{\pi(t)}=[{\|\pi_{0}\|}_{W_{p}^1(\Omega)}^p +
cR_{(\overline{\sigma},t)} (R_{(\overline{T},t)}t^{\frac{q-1}{q}}+
R_{(\overline{\pi},t)}t)]\times
\end{equation}
$$
\times
\exp\big(c(R_{(\overline{v},t)}t^{\frac{p-1}{p}}+R_{(\overline{\sigma},t)}
(R_{(\overline{T},t)}t^{\frac{q-1}{q}}+R_{(\overline{\pi},t)}t))\big),
$$
$$
R_{(\overline T, t)} = \| \overline T \|_{W^{2,1}_q(Q_t ) } ,
\qquad
R_{(\overline \pi , t)} = \| \overline \pi \|_{C^0 ([0, t ] ;
W^1_p( \Omega ) )},
$$
$$
 \qquad R_{( {\overline \sigma  ,t} )}  = \|
\overline \sigma \|_{C^0 ( [0,t],W_p^1 ( D_2))}.
$$
\end{lem}

{\bf Proof.}  See \cite{[SF]}. $\square$

\begin{lem}\label{lem-sigma}
Let be $\overline v \in \Theta^{(v)}_{t_1} $, $\overline T \in
\Theta^{(T)}_{t_1}$, $ \overline \pi \in C^0([0,t_1];
W^1_p(\Omega)) $ , $ \overline \sigma \in C^0 ([0, t_1 ];
W^1_p(\mathbf{R}_+ \times \Omega ) ) $. Moreover we assume
that $
\partial_m \overline \sigma \in C^0 ([0, t_1 ];
L^\infty(\mathbf{R}_+ \times \Omega ) ) $. Then there exists a
positive constant $\overline M _1$ such that, provided that $\overline
\sigma(m,\cdot,\cdot)=0$ for $m\notin \, \big] {\frac{{\overline m _a }}{2} ,
\overline M _1 } \big[ $, the equation
\begin{equation}\label{ex(5.16)}
\frac{\partial \sigma }{\partial t} + \nabla_{(m,x)} \cdot (\sigma
\widetilde U_{4l} (\overline u,\overline T,\overline{\pi} ) ) =
\end{equation}
$$
= \big[ h_{gl} ( \overline T ,\overline{\pi} ;m ) + B_1 (
\overline{\sigma} ; m)
 - g_1 ( m )
[ \overline{\pi} - \overline{ \pi}_{vs(l)}( \overline T ) ]^- \big]
\sigma +
$$
$$
+ g_0 ( m ) [ N^*  - \widetilde N ( \overline{\sigma}  ) ]^+ [
\overline{\pi} - \overline{ \pi}_{vs(l)}( \overline T )]^+ + B_2 (
\overline {\sigma} ; m ),
$$
with the initial condition
\eqref{cond-sigma}--\eqref{4.15-sigma-nu-bis} admit a unique
solution $\sigma \in C^0 ([0, t_1 ]; W^1_p(\mathbf{R}_+ \times
\Omega ) ) $ satisfying the following relations
\begin{equation}\label{ex(new-5.18)}
{  \partial_m \sigma \in C^0 ([0, t_1 ]; L^\infty(\mathbf{R}_+
\times \Omega ) ) ,  }
\end{equation}
\begin{equation}\label{ex(5.18)}
\sigma( \cdot,m,\cdot)  =0 \quad  \ \mbox{for} \ \ m \notin \, \big]
{\frac{{\overline m _a }}{2} ,\overline M _1 } \big[,
\end{equation}
\begin{equation}\label{ex(5.19)}
\| \sigma (t,\cdot , \cdot ) \|_{W^1_p ( \mathbf{R}_+ \times
\Omega )}^{\, p} = \| \sigma (t,\cdot , \cdot ) \|_{W^1_p (
\Omega_{\overline M_1})}^{\, p} \leq q_\sigma (t) ,
\end{equation}
\begin{equation}\label{ex(5.19-deltam-sigma)}
\left\| {\partial _m \sigma ( t,\cdot )} \right\|_{L^\infty
(\mathbf{R}_+ \times \Omega )}  \le \Big [\left\| {\partial _m
\sigma _0 } \right\|_{L^\infty (\mathbf{R}_+ \times \Omega )}  +
c\int\limits_0^t {\Big(1 + \left\| {\overline \pi  (s, \cdot )}
\right\|_{L^\infty  (\Omega )}^2    }
+
\end{equation}
$$
+\left\| {\overline \sigma (s,
\cdot )} \right\|_{L^\infty  (\mathbf{R}_+ \times \Omega)}^2+
\left\| {\sigma (s, \cdot )} \right\|_{L^\infty (\mathbf{R}_+
\times \Omega )}^2  + \left\| {\partial _m \overline \sigma  (
s,\cdot )} \right\|_{L^\infty  (\mathbf{R}_+ \times \Omega)}^2
\Big)ds\Big] \times
$$
$$
\times \exp \Big [c\int\limits_0^t {\Big (1 + \left\| {\nabla _x
\cdot \overline v (s, \cdot )} \right\|_{L^\infty  (\Omega )}  +
\left\| {\overline \pi  (s, \cdot )} \right\|_{L^\infty  (\Omega
)} + \left\| {\overline \sigma  (s, \cdot )} \right\|_{L^\infty
(\mathbf{R}_+ \times \Omega)} } \Big )ds \Big],
$$
where
\begin{equation}\label{ex(5.20)}
q_\sigma (t)= \big\{ \| \sigma_0 \|_{W_p^1 (\Omega_{\overline M_1})}^p + c \big[ (1 +
R_{( \overline \pi , t )}^2 + R_{( \overline \sigma , t )}^2)t
 +
R_{ ( \overline v ,t )} t^{\frac{{p - 1}}{p}} + R_{ ( \overline T
,t )}^2 t^{\frac{{q - 2}}{q}} \big] \big\} \times
\end{equation}
$$
\times \exp \big\{ c \big[ (1 + R_{( \overline \pi , t )}^2 + R_{(
\overline \sigma , t )}^2 + R_{ ( \overline v ,t )} t^{\frac{{p -
1}}{p}} + R_{ ( \overline T ,t )}^2 t^{\frac{{q - 2}}{q}} \big]
\big\} .
$$

\end{lem}

{\bf Proof.}\ For the proof of the existence and uniqueness of the
solution in $ C^0([0,t_1];W^1_p(\Omega))$ and the relation
\eqref{ex(5.18)}, see \cite{[SF]}. The relation \eqref{ex(new-5.18)}
will follow from \eqref{ex(5.19-deltam-sigma)}.

For the inequality \eqref{ex(5.19)},
by taking into account  \eqref{4.15-geo-a}, \eqref{u-bar-v-bar}
and \eqref{ex(5.18)}, from \eqref{eq-sigma-lin} we deduce
\begin{equation}\label{5.45-gauss-1}
\int_{\Omega_{\overline M_1} } {\sigma^{p - 1} \nabla_{( m,x)} } \sigma \cdot
\widetilde U_{4l} ( \overline u ,\overline T ,\overline \pi )dmdx
=- \frac{1}{p}\int_{\Omega_{\overline M_1} } \sigma^p \nabla_{ ( m,x )} \cdot
\widetilde U_{4l} ( \overline u ,\overline T ,\overline \pi )dmdx
,
\end{equation}
\begin{equation}\label{5.45-gauss-2}
\int_{\Omega_{\overline M_1} } \left| \nabla_{( m,x)} \sigma \right|^{p - 2}
\nabla_{(m,x)} \sigma  \cdot \big( {\widetilde U_{4l} ( \overline
u ,\overline T ,\overline \pi ) \cdot \nabla_{(m,x)} } \big)
\nabla_{(m,x)} \sigma dmdx =
\end{equation}
$$
=  - \frac{1}{p}\int_{\Omega_{\overline M_1} } \left| \nabla_{\left( {m,x} \right)}
\sigma \right|^p \nabla_{( m,x)}  \cdot \widetilde U_{4l} (
\overline u ,\overline T ,\overline \pi   )dmdx .
$$

Multiplying the equation \eqref{ex(5.16)} by $\sigma^{p - 1} $ and
integrating it on ${\Omega_{\overline M_1} }$, thanks to
\eqref{5.45-gauss-1}, we obtain
\begin{equation}\label{5.50-norma-sigma}
\frac{d}{{dt}} \| \sigma \|_{L^p ( \Omega_{\overline M_1} )}^p  = ( 1 - p )\int_{\Omega_{\overline M_1}
} \sigma^p \nabla _{(m,x)} \cdot \widetilde U_{4l} ( \overline u
,\overline T , \overline \pi ) dmdx  +
\end{equation}
$$
+ p\int_{\Omega_{\overline M_1} } \sigma^p a^*_1 (t,m,x) dmdx + p\int_{\Omega_{\overline M_1} }
\sigma^{p - 1} b_1 (t,m,x ) dmdx
$$
where
$$
a^*_1 ( t,m,x ) = h_{gl} (\overline T ,\overline{\pi} ;m) + B_1
(\overline{\sigma} ; m )
 - g_1 ( m )
[ \overline{\pi}  - \overline{ \pi}_{vs}(\overline T ) ]^- ,
$$
$$
b_1 ( {t,m,x} ) =
g_0 ( m ) [ N^*  - \widetilde N ( \overline{\sigma} ) ]^+
[ \overline{\pi}  - \bar \pi_{vs(l)}( \overline T )]^+
+ B_2 ( \overline{\sigma} ; m ).
$$
On the other hand, applying the differential operator $ |
\nabla_{( m,x )} \sigma |^{p - 2} \nabla _{(m,x)} \sigma \cdot
\nabla _{(m,x)} $ to the equation \eqref{ex(5.16)} and integrating
it on ${\Omega_{\overline M_1} }$, we have
\begin{equation}\label{5.50-norma-grad-sigma}
\frac{d}{{dt}} \| \nabla_{(m,x)} \sigma \|_{L^p (\Omega_{\overline M_1})}^p  = -
p\int_{\Omega_{\overline M_1} } { | \nabla_{( m,x)} \sigma |^{p - 2} \nabla_{( m,x)}
\sigma  \cdot \nabla_{( m,x)} \big[ \nabla_{( m,x)} \cdot ( \sigma
\widetilde U_{4l} ( \overline u ,\overline T , \overline \pi ) )
\big]} dmdx +
\end{equation}
$$
+ p\int_{\Omega_{\overline M_1} } | \nabla_{( m,x)} \sigma |^{p - 2} \nabla_{( m,x)}
\sigma  \cdot \nabla_{( m,x)} [ a^*_1 (m,x,t )\sigma + b_1 ( t,m,x
) ] dmdx .
$$
Remembering the definitions of $ a_1^* (t,m,x )$ and $b_1 ( t,m,x
) $ and using repeatedly  the Sobolev and H\"older inequalities,
we deduce
\begin{equation}\label{5.60-norma-1}
\Big|\int_{\Omega_{\overline M_1} } \sigma^p \nabla_{(m,x)} \cdot \widetilde U_{4l} (
\overline u , \overline T ,\overline \pi ) dmdx\Big| + \Big|
\int_{\Omega_{\overline M_1} } \sigma^p a^*_1 ( {t,m,x} ) dmdx \Big| \leq
\end{equation}
$$
\leq c \big( 1 + \| \overline \pi \|_{W_p^1 ( \Omega )} + \|
\overline T \|_{W_q^2 ( \Omega )} + \| \overline u \|_{W_p^2 (
\Omega_{\overline M_1})} + \| \overline \sigma \|_{W_p^1 ( \Omega_{\overline M_1})}
 \big)
\| \sigma \|_{L^p ( \Omega_{\overline M_1} )}^p ,
$$
\begin{equation}\label{5.60-norma-3}
\Big| \int_{\Omega_{\overline M_1} } \sigma^{p - 1} b_1 ( {t,m,x} )dmdx \Big| \leq \|
\sigma \|_{L^p ( \Omega_{\overline M_1} )}^{p-1} \| b_1 \|_{L^p ( \Omega_{\overline M_1} )} \leq
\end{equation}
$$
\leq c \big( 1+ \| \overline \pi \|_{L^p ( \Omega )}^2 +\|
\overline T \|_{W_q^1 ( \Omega )}^2 + \| \overline \sigma
\|_{W^1_p ( \Omega_{\overline M_1})}^2 \big) (\| \sigma \|_{L^p ( \Omega_{\overline M_1} )}^p + 1 ) ,
$$
\begin{equation}\label{5.70-norma-1}
\Big| \int_{\Omega_{\overline M_1} } { | \nabla_{( m,x)} \sigma |^{p - 2} \nabla_{(
m,x)} \sigma  \cdot \nabla_{( m,x)} \big[ \nabla_{( m,x)} \cdot (
\sigma \widetilde U_{4l} ( \overline u ,\overline T , \overline
\pi ) ) \big]} dmdx \Big| \leq
\end{equation}
$$
\leq c \big( 1 + \| \overline \pi  \|_{W_p^1 ( \Omega )} + \|
\overline T \|_{W_q^2 ( \Omega )} + \| \overline u \|_{W_p^2 (
\Omega_{\overline M_1}) } \big) \| \sigma \|_{W_p^1 (\Omega_{\overline M_1})}^p ,
$$
\begin{equation}\label{5.70-norma-2}
\Big| \int_{\Omega_{\overline M_1} } | \nabla_{( m,x)} \sigma |^{p - 2} \nabla_{ (
m,x )} \sigma \cdot \nabla_{ ( m,x )} ( a^*_1 ( t,m,x )\sigma
)dmdx \Big| \leq
\end{equation}
$$
\leq c (1 + \| \overline \pi \|_{W_p^1 ( \Omega )} + \| \overline
T \|_{W_q^2 ( \Omega )} + \| \overline \sigma \|_{W_p^1 ( \Omega_{\overline M_1} )}
 )
\| \sigma \|_{W_p^1 ( \Omega_{\overline M_1} )}^p ,
$$
\begin{equation}\label{5.70-norma-3}
\Big| \int_{\Omega_{\overline M_1} } | \nabla_{ ( m,x )} \sigma |^{p - 2} \nabla_{ (
m,x )} \sigma \cdot \nabla_{( m,x )} b_1 ( t,m,x )dmdx \Big| \leq
\end{equation}
$$
\leq c (1 + \| \overline \pi \|_{W_p^1 ( \Omega )}^2 + \|
\overline T \|_{W_q^2 ( \Omega )}^2 + \| \overline \sigma
\|_{W_p^1 ( \Omega_{\overline M_1} )}^2
 )
(\| \sigma \|_{W^1_p ( \Omega_{\overline M_1} )}^p + 1 ) .
$$

From \eqref{5.50-norma-sigma}--\eqref{5.70-norma-3}, taking into
account \eqref{u-bar-v-bar}, we obtain
$$
\frac{d}{{dt}} \| \sigma \|_{W_p^1 ( \Omega_{\overline M_1} )}^p \leq c (1 + \|
\overline v \|_{W_p^2 ( \Omega )} + \| \overline \pi \|_{W_p^1 (
\Omega )}^2 + \| \overline T \|_{W_q^2 ( \Omega )}^2 + \|
\overline \sigma \|_{W_p^1 ( \Omega_{\overline M_1} )}^2
 )
(\| \sigma \|_{W^1_p ( \Omega_{\overline M_1} )}^p + 1 ) ,
$$
from this inequality we deduce \eqref{ex(5.19)} with $q_\sigma (t)$ defined as
in \eqref{ex(5.20)}.

To prove the inequality \eqref{ex(5.19-deltam-sigma)}, we consider $r
\geq p$. Applying the operator $\left| {\partial _m \sigma }
\right|^{r - 2} \partial _m \sigma \partial _m$ to the equation
\eqref{ex(5.16)} and integrating it on $\Omega_{\overline M_1}$,
we have
\begin{equation}\label{5.70-norma-Lr}
\frac{d}{{dt}}\left\| {\partial _m \sigma } \right\|_{L^r (\Omega_{\overline M_1})}
\le c\Big [\Big (1 + \left\| {\nabla _x  \cdot \overline u (s, \cdot)} \right
\|_{L^\infty  (\Omega_{\overline M_1})}  + \left\| {\overline \pi  }
\right\|_{L^\infty  (\Omega )}
 + \left\| {\overline \sigma  } \right\|_{L^\infty  (\Omega_{\overline M_1})} )
 \left\| {\partial _m \sigma } \right\|_{L^r (\Omega_{\overline M_1})}  +
\end{equation}
$$
 + \left\| {\overline \sigma  } \right\|_{L^\infty  (\Omega_{\overline M_1} )}
\Big(\left\| {\overline \sigma  } \right\|_{L^r (\Omega_{\overline M_1} )} + \left\|
{\partial _m \overline \sigma  } \right\|_{L^r (\Omega_{\overline M_1} )} \Big) +
\Big(1 + \left\| {\overline \pi  } \right\|_{L^\infty  (\Omega )}
+ \left\| {\overline \sigma  } \right\|_{L^\infty  (\Omega_{\overline M_1} )}
\Big)\left\| \sigma  \right\|_{L^r (\Omega_{\overline M_1} )} \Big].
$$
Applying Gronwall's lemma and taking the limit for $r \to \infty$
we obtain \eqref{ex(5.19-deltam-sigma)}. \ $\square$

\section{Equations for the water densities with given velocities and temperature}

To prove the existence and uniqueness of the solution
$(\pi,\sigma)$ to the nonlinear equation system
\eqref{eq-pi-lin}--\eqref{eq-sigma-lin} and obtain its estimates
with given temperature $T=\overline T$ and velocity $v = \overline
v$, we use the following lemma.

\begin{lem}\label{lem-non-lin-1}
Let be $ \overline v \in \Theta^{(v)}_{t_1} $, $ \overline T \in
\Theta^{(T)}_{t_1} $ and $R_0>0$. We assume that
\begin{equation}\label{ex(6.1)-borne}
{\|\overline{v}\|}_{{W_{p}^{2,1}}(Q_{t_1})} \leq
R_0,\qquad{\|\overline{T}\|}_ {{W_{q}^{2,1}}(Q_{t_1})} \leq R_0.
\end{equation}
Then there exists $ t_2 = t_2(R_0)$, $0< t_2 \leq t_1$, such that,
if $\overline \pi \in C^0([0, t_2]; W^1_p( \Omega ) ) ,$
$\overline \sigma \in C^0([0, t_2]; W^1_p( \mathbf{R}_+ \times
\Omega ) )$ and $ \partial_m \overline \sigma \in C^0 ([0, t_1 ];
L^\infty(\mathbf{R}_+ \times \Omega ) )$ with the conditions
\begin{equation}\label{ex(6.2)} \| \overline \pi
\|_{C^0([0, t_2]; W^1_p( \Omega ) )} \leq \| \pi_0 \|_{ W^1_p(
\Omega ) } + 1 , \qquad \| \overline \sigma \|_{C^0([0, t_2];
W^1_p( \mathbf{R}_+ \times \Omega ))} \leq \| \sigma_0 \|_{ W^1_p(
\mathbf{R}_+ \times \Omega )} + 1 ,
\end{equation}
$$
 \|  \partial_m \overline \sigma \|_{C^0([0, t_2];
L^\infty( \mathbf{R}_+ \times \Omega ))} \leq \| \partial_m
\sigma_0 \|_{ L^\infty( \mathbf{R}_+ \times \Omega )} + 1 ,
$$
$$
\overline \sigma(\cdot,m,\cdot) = 0 \quad \mbox{for} \ \ m \not \in \;
\big] \frac{{\overline m _a }}{2}, \overline M_1 \big[ ,
$$
then the solution $(\pi , \sigma  )$ of the equations
\eqref{ex(5-eq-pi-bar2)}, \eqref{ex(5.16)}  with the initial
conditions \eqref{cond-pi}--\eqref{cond-sigma} satisfies the
conditions
\begin{equation}\label{ex(6.3)-stime}
\| \pi \|_{C^0([0, t_2]; W^1_p( \Omega ) )} \leq \| \pi_0 \|_{
W^1_p( \Omega ) } + 1 , \qquad \| \sigma \|_{C^0([0, t_2]; W^1_p(
\mathbf{R}_+ \times \Omega ))} \leq \| \sigma_0 \|_{ W^1_p(
\mathbf{R}_+ \times \Omega )} + 1 ,
\end{equation}
$$
 \|  \partial_m \sigma \|_{C^0([0, t_2];
L^\infty( \mathbf{R}_+ \times \Omega ))} \leq \| \partial_m
\sigma_0 \|_{ L^\infty( \mathbf{R}_+ \times \Omega )} + 1 ,
$$
$$
\sigma(\cdot,m,\cdot) = 0 \quad \mbox{for} \ \ m \not \in \;
\big]\frac{{\overline m _a }}{2}, \overline M_1 \big[ .
$$
\end{lem}
{\bf Proof.} \ The lemma follows from the relations
\eqref{ex(5.14)}, \eqref{ex(5.18)}, \eqref{ex(5.19)}, \eqref{ex(5.19-deltam-sigma)}
(see also \eqref{ex(5.15)}, \eqref{ex(5.20)}). \ $\square$


\begin{lem}\label{lem-non-lin-2}
Let be $\overline v$, $\overline T$, $\overline M_1$, $R_0$, $t_2
= t_2(R_0)$ as in Lemma \ref{lem-non-lin-1}. Then there exists
$t_3\in \, ]0, t_2]$ such that the equation system
\eqref{eq-pi-lin}--\eqref{eq-sigma-lin} with the initial conditions
\eqref{cond-pi}--\eqref{cond-sigma} admit a unique solution $(\pi
, \sigma)$ $\in$  $C^0 ([0, t_3 ]; W^1_p(\Omega ) )$ $\times$ $
C^0 ([0, t_3 ]; W^1_p(\mathbf{R}_+ \times \Omega ) )) $,
satisfying the conditions
\begin{equation}\label{ex(6.10)}
\| \pi \|_{C^0([0, t_3]; W^1_p( \Omega ) )} \leq \| \pi_0 \|_{
W^1_p( \Omega ) } + 1 , \qquad \| \sigma \|_{C^0([0, t_3]; W^1_p(
\mathbf{R}_+ \times \Omega ))} \leq \| \sigma_0 \|_{ W^1_p(
\mathbf{R}_+ \times \Omega )} + 1 ,
\end{equation}
$$
 \|  \partial_m \sigma \|_{C^0([0, t_2];
L^\infty( \mathbf{R}_+ \times \Omega ))} \leq \| \partial_m
\sigma_0 \|_{ L^\infty( \mathbf{R}_+ \times \Omega )} + 1 ,
$$
$$
\sigma(\cdot,m,\cdot) = 0 \quad \mbox{for} \ \ m \not \in \;
\big]\frac{{\overline m _a }}{2}, \overline M_1 \big[ .
$$

\end{lem}

{\bf Proof.} \ For $0< t \leq t_2$, we define the set
\begin{equation}\label{ex(6.5)}
A_{[ t ]}  = \left\{ {( {\pi ,\sigma  } )} \ \mbox{satisfies the
conditions $(C.1)$--$(C.4)$}  \right\} ,
\end{equation}
where
$$
\pi \in C^0 ([0, t ]; W^1_p(\Omega ) ), \quad \sigma\in C^0 ([0, t
]; W^1_p(\mathbf{R}_+ \times \Omega ) ), \leqno (C.1)
$$
$$
\pi(0,\cdot )=\pi_0 (\cdot) , \quad \sigma (0,\cdot,\cdot)
=\sigma_0(\cdot,\cdot )  , \leqno (C.2)
$$
$$
\| \pi \|_{C^0([0, t]; W^1_p( \Omega ) )} \leq \| \pi_0 \|_{
W^1_p( \Omega ) } + 1 , \qquad \| \sigma \|_{C^0([0, t]; W^1_p(
\mathbf{R}_+ \times \Omega ))} \leq \| \sigma_0 \|_{ W^1_p(
\mathbf{R}_+ \times \Omega )} + 1 , \leqno (C.3)
$$
$$ \|  \partial_m \sigma \|_{C^0([0, t_2];
L^\infty( \mathbf{R}_+ \times \Omega ))} \leq \| \partial_m
\sigma_0 \|_{ L^\infty( \mathbf{R}_+ \times \Omega )} + 1 ,
$$
$$
\sigma(t',m,x) = 0 \quad \mbox{for} \ \ m \not \in \;
\big]\frac{{\overline m _a }}{2}, \overline M_1 \big[ , \ x \in \Omega,
\ 0\leq t' \leq t . \leqno (C.4)
$$

From Lemma 6.1, the map $ G_{1,t} : A_{[ t ]} \to A_{[ t ]} $ that
sends  $( \overline \pi ,\overline \sigma ) \in A_{[ t ]}$ to the
solution $(\pi , \sigma ) = G_{1,t} (
\overline \pi , \overline \sigma   )$ of the equations
\eqref{ex(5-eq-pi-bar2)}, \eqref{ex(5.16)} is well defined.

Now, we prove that there exists  $t_3$, $0 < t_3 \leq t_2$, such
that $G_{1,t_3}$ is a contraction with respect to the metric of
$C^0 ( [0,t_3];L^p ( \Omega)) \times C^0 ( [0,t_3];L^p
(\Omega_{\overline M_1}))  $ on the closed set $A_{[ t_3]}$. For
this, we consider two elements $ (\overline \pi_1 , \overline
\sigma_1 )$, $ (\overline \pi_2 , \overline \sigma_2 )$ of $A_{[
t]}$ with their values by the map $G_{1,t}$ are $( \pi_i ,\sigma_i
) = G_{1,t} ( \overline \pi_i ,\overline \sigma_i )$, $i = 1,2$;
we denote by
\begin{equation}\label{Sigma-Pi-Xi}
\overline \Pi  = \overline \pi_2 - \overline \pi_1 , \quad
\overline \Sigma  = \overline \sigma_2 - \overline \sigma_1 ,\quad
\Pi  = \pi_2  - \pi_1 , \quad \Sigma  = \sigma_2 - \sigma_1 .
\end{equation}
The difference between the equations \eqref{ex(5-eq-pi-bar2)},
\eqref{ex(5.16)} with $ (\overline \pi_2 , \overline \sigma_2  )$
and $ (\overline \pi_1 , \overline \sigma_1  )$ gives
\begin{equation}\label{6.5-pi-maius}
\frac{{\partial \Pi }}{{\partial t}} + \nabla  \cdot \left( {\Pi
\overline v } \right) = H_{gl} \left( {\overline T , \overline
\pi_{1 } ,\overline \sigma_1 } \right) - H_{gl} \left( {\overline
T ,\overline \pi_{2} , \overline \sigma_2 } \right) ,
\end{equation}
\begin{equation}\label{6.5-sigma-maius}
\frac{ \partial \Sigma }{ \partial t} + \frac{\partial }{{\partial
m}} \big[ mh_{gl} ( \overline T ,\overline \pi_{2 } ;m ) \Sigma
\big] + \nabla  \cdot ( \Sigma \overline u ) =
\end{equation}
$$
= \big[ h_{gl} ( \overline T ,\overline \pi_{2 } ;m ) + B_1 (
\overline \sigma_2 ; m )  - g_1 ( m ) [ \overline \pi_2  - \overline
{\pi}_{vs} ( \overline T )]^- \big]\Sigma +
$$
$$
+ \frac{\partial }{{\partial m}}\big\{ m [ h_{gl} ( \overline T
,\overline \pi_{1 } ;m ) - h_{gl} ( \overline T ,\overline \pi_{2}
;m )] \sigma_1 \big\} + [ h_{gl} ( \overline T ,\overline \pi_{2 }
;m ) - h_{gl} ( \overline T ,\overline \pi_{1 } ;m ) ] \sigma_1 +
$$
$$
+ [ B_1 ( \overline \sigma_2 ;m ) - B_1 ( \overline \sigma_1 ;m) ]
\sigma_1  +
$$
$$
+ g_1 ( m )\big( [ \overline \pi_1  - \overline{ \pi}_{vs} ( \overline
T ) ]^- - [ \overline \pi_2  - \overline{\pi}_{vs} ( \overline T ) ]^
-  \big) \sigma _1 +
$$
$$
+ g_0 ( m )\big( [ N^*  - \widetilde N ( \overline \sigma_2 ) ]^+
[ \overline \pi_2 - \overline{ \pi}_{vs} ( \overline T )]^+ - [ N^*  -
\widetilde N ( \overline \sigma_1)]^+ [ \overline \pi_1 - \overline
{\pi}_{vs} ( \overline T )]^+ \big) +
$$
$$
+ B_2 ( \overline \sigma_2 ;m) - B_2 ( \overline \sigma_1 ;m) .
$$
Multiplying the equations
\eqref{6.5-pi-maius}--\eqref{6.5-sigma-maius} respectively by $|
\Pi |^{p - 1} $ and $| \Sigma |^{p - 1} $, integrating the first
on $\Omega $ and the second on  $\Omega_{\overline M_1}$ and taking
into account the conditions \eqref{ex(6.2)} and the estimates
already known, we deduce the following inequalities
\begin{equation}\label{eq-norma-Pi}
\frac{d}{dt} \| \Pi \|_{L^p ( \Omega )}^p \leq c \| \overline v
\|_{W^2_p (\Omega)} \| \Pi \|_{L^p ( \Omega )}^{p} +
\end{equation}
$$
+ c (1 + \| \overline T \|_{W^1_q (\Omega)}) \big( \| \overline
\Pi \|_{L^p ( \Omega )} + \| \overline \Sigma \|_{L^p ( \Omega_{\overline M_1} )}
\big) \| \Pi \|_{L^p ( \Omega )}^{p-1} ,
$$
\begin{equation}\label{eq-norma-Sigma}
\frac{d}{dt} \| \Sigma \|_{L^p ( \Omega_{\overline M_1} )}^p \leq c (1 + \| \overline
T \|_{W^1_q (\Omega)} + \| \overline v \|_{W^2_p (\Omega)} ) \|
\Sigma \|_{L^p ( \Omega_{\overline M_1} )}^{p} +
\end{equation}
$$
+ c (1 + \| \overline T \|_{W^1_q (\Omega)}) \big( \| \overline
\Pi \|_{L^p ( \Omega )} + \| \overline \Sigma \|_{L^p ( \Omega_{\overline M_1} )}
\big) \| \Sigma \|_{L^p ( \Omega_{\overline M_1} )}^{p-1} .
$$
It is useful to recall that $\|\partial_m \sigma_0 \|_{ L^\infty( \mathbf{R}_+ \times
\Omega )} < \infty$, the inequality \eqref{eq-norma-Sigma} is
obtained without introducing a regularization of $\pi$. Indeed,
the term that had requested a regularization of $\pi$ in \cite{[SF]} can
be treated in the following way
$$
\Big|\int\limits_{\Omega_{\overline M_1} } {\left| \Sigma  \right|^{p - 1}
\partial _m \Big\{m \Big[h_{gl} (T,\overline \pi  _1 ;m) - h_{gl}
(T,\overline \pi  _2;m )\Big]\sigma _1 \Big\}dmdx \Big|} \le
$$
$$
\le c\int\limits_{\Omega_{\overline M_1} } {\left| \Sigma  \right|^{p - 1} \left|
{\overline \Pi  } \right|(\left| {\sigma _1 } \right| + \left|
{\partial _m \sigma _1 } \right|)dmdx}  \le
$$
$$
\le c \Big(\left\| {\sigma _1 (t, \cdot )} \right\|_{L^\infty
(\Omega_{\overline M_1} )} + \left\| {\partial _m \sigma _1 (t, \cdot )}
\right\|_{L^\infty (\Omega_{\overline M_1} )} \Big) \left\| {\Sigma (t, \cdot )}
\right\|_{L^p (\Omega_{\overline M_1} )}^{p-1} \left\| {\overline \Pi  (t, \cdot )}
\right\|_{L^p (\Omega )}.
$$

Multiplying
\eqref{eq-norma-Pi}--\eqref{eq-norma-Sigma}  respectively by $\| \Pi \|_{L^p (
\Omega )}^{1-p} $ and $\| \Sigma \|_{L^p ( \Omega )}^{1-p} $ and
taking into account \eqref{ex(6.1)-borne}, we deduce
\begin{equation}\label{6.10-pi-norma}
\| \Pi (t, \cdot )\|_{L^p ( \Omega )} + \| \Sigma (t, \cdot)\|_{L^p ( \Omega_{\overline M_1} )}
 \leq
\end{equation}
$$
\leq c e^{c( t+ t^{\frac{q-1}{q}}+ t^{\frac{p-1}{p}})} ( t+
t^{\frac{q-1}{q}}) \big( \| \overline \Pi  \|_{C^0 ( [0,t];L^p (
\Omega ))} + \| \overline \Sigma  \|_{C^0 ( [0,t];L^p ( \Omega_{\overline M_1} ))}
 \big).
$$
It is clear that there exists $t_3\in \, ]0, t_2]$ such that
$$
c e^{c( t_3+ t_3^{\frac{q-1}{q}}+ t_3^{\frac{p-1}{p}} )} ( t_3 +
t_3^{\frac{q-1}{q}}) < \frac{1}{2} .
$$
Therefore, remembering the definition of $G_{1,t}$ and
\eqref{Sigma-Pi-Xi}, we deduce that the map $G_{1,t_3}$ restricted
to $A_{[ t_3 ]} $ is a contraction with respect to the metric of
$C^0 ([ 0,t_3];L^p ( \Omega )) \times C^0 ([ 0,t_3];L^p (
\Omega_{\overline M_1}))$. \ $\square$

\medskip

\begin{lem}
Let be $ \overline v \in \Theta^{(v)}_{t_1} $, $ \overline T \in
\Theta^{(T)}_{t_1} $ verifying \eqref{ex(6.1)-borne} and
$\varrho$, $\pi$, $\sigma$  the solutions of the equations
\eqref{eq-varrho-lin}--\eqref{eq-sigma-lin} with the initial
conditions \eqref{cond-varrho}--\eqref{cond-sigma}. Under the same
hypotheses of Lemma 6.2, there exists $t_4 \in
\, ]0,t_3]$ such that for $ t \in [0, t_4]$
we have
\begin{equation}\label{ex(6.40)a}
\| \varrho (t,\cdot )\|_{W^1_p }^p  \leq 2 \| \varrho_0 (\cdot
)\|_{W^1_p }^p , \qquad {1\over 2}  \inf_{x' \in
\Omega}\varrho_0(x') \leq \varrho (t,x) \leq  2\sup_{x' \in
\Omega}\varrho_0(x') \quad x \in \Omega,
\end{equation}
\begin{equation}\label{ex(6.40)b}
0  \leq \pi (t,x) \leq  \sup_{x' \in \Omega}\pi_0(x') + 1\quad x \in \Omega.
\end{equation}
\end{lem}

{\bf Proof.}
See \cite{[SF]}.
\ $\square$

\section{ Linear equation for the velocity and temperature}

In this section we study the linearized equations of
\eqref{ex(2.1)} and \eqref{ex(2.4)}. We assume that $\overline v\in
\Theta^{(v)}_{t_1}$ and $\overline T \in \Theta^{(T)}_{t_1}$ and
consider the linear equations in $v$ and $T$
\begin{equation}\label{7-eq-v}
(\varrho + \pi )  {\partial v \over \partial t } -\eta \Delta v -
\big( \zeta + {\eta \over 3 } \big) \nabla ( \nabla \cdot v ) =
\end{equation}
$$
= -(\varrho + \pi ) ( \overline v \cdot \nabla ) \overline v - R_0
\nabla ( ( {\varrho \over \mu_a} + {\pi \over \mu_h} ) \overline T
) - \Big[ \int_0^\infty ( \sigma + \nu )dm + \varrho  + \pi \Big]
\nabla \Phi ,
$$
\begin{equation}\label{7-eq-T}
(\varrho +\pi)c_v  {\partial T \over \partial t } -\kappa \Delta T
=-(\varrho +\pi)c_v \sum_{j=1}^3 \overline v_j {\partial \overline
T \over \partial x_j } - R_0 ( {\varrho \over \mu_a} + {\pi \over
\mu_h} ) \overline T \nabla \cdot \overline v +
\end{equation}
$$
+ \eta \sum_{i,j=1}^3 \Big( {\partial \overline v_i \over \partial
x_j} + {\partial \overline v_j \over \partial x_i} - {2 \over 3}
\delta_{ij} \nabla \cdot \overline v \Big) {\partial \overline v_i
\over \partial x_j} + \zeta (\nabla \cdot \overline v)^2 - \nabla
\cdot \mathcal{E} +
$$
$$
+ L_{gl} H_{gl} (\overline T, \pi , \sigma ) + L_{ls}
H_{ls}(\overline T, \sigma , \nu ) + L_{gs} H_{gs}(\overline T,
\pi , \nu ) ,
$$
where $\varrho $, $\pi $, $\sigma$ are the solutions of the equations
\eqref{eq-varrho-lin}--\eqref{eq-sigma-lin} with the initial
conditions \eqref{cond-varrho}--\eqref{cond-sigma}; their
existence and uniqueness are proved in Lemma 6.1 and Lemma 6.2.

As in \cite{[So76]}, we introduce the following auxiliary
functions
\begin{equation}\label{ex(7.4)}
V_{(p,v)}(t) = \| v \|_{ W^{2,1}_p(Q_{t})}^p + \sup_{0\leq t' \leq
t} \| v( t',\cdot) \|_{ W^{2-{2\over p}}_p(\Omega)}^p ,
\end{equation}
\begin{equation}\label{ex(7.6)}
V_{(q,T)}(t) = \| T \|_{ W^{2,1}_q(Q_{t})}^q + \sup_{0\leq t' \leq
t} \| T(t',\cdot) \|_{ W^{2-{2\over q}}_q(\Omega)}^q
\end{equation}
(for the general theory about the concerned functional spaces, see
also \cite{[ISo]}, \cite{[So64]}, \cite{[So65]}).

\begin{lem}
Let be $\overline v\in \Theta^{(v)}_{t_1}$, $\overline T \in
\Theta^{(T)}_{t_1}$and $( \varrho ,\pi ,\sigma)$ the solution of
the equation system \eqref{eq-varrho-lin}--\eqref{eq-sigma-lin}
with the initial conditions \eqref{cond-varrho}--\eqref{cond-sigma}
given in Lemma 5.1 and Lemma 6.2. Then the equations
\eqref{7-eq-v} and \eqref{7-eq-T}, with the conditions
\eqref{cond-v}--\eqref{cond-T}, admit a unique solution
\begin{equation}\label{ex(7.7)}
v \in W^{2,1}_p(Q_{ t_4}) , \qquad T \in W^{2,1}_q(Q_{ t_4}).
\end{equation}
Moreover we have
\begin{equation}\label{ex(7.8)}
V_{(p,v)}(t) \leq c \Big( \| v_0 \|_{ W^{2-{2\over
p}}_p(\Omega)}^p + \int_0^t (1 + V_{(p,\overline v)}(t')^2 ) dt'+
t^{\frac{2q-p}{q}}V_{(q,\overline T )}(t) \Big) ,
\end{equation}
\begin{equation}\label{ex(7.10)}
V_{(q,T)}(t) \leq  c \Big[ \| T_0\|_{ W^{2-{2\over
q}}_q(\Omega)}^q   +
t^{ \frac{2(p-q)}{ p} } V_{(p,\overline v )}(t)^{\frac{2q}{p}} +
\end{equation}
$$
+ \int_0^t (1 + V_{(q,\overline T )}(t') V_{(p,\overline v
)}(t')^{q / p} + \| \nabla \cdot \mathcal{E}\|^q_{L^q(\Omega)} )
dt' \Big]
$$
for $ 0<t \leq t_4 $ ($t_4$ is defined in Lemma 6.3).
\end{lem}

{\bf Proof.} \
According to Theorem 9.1 of the  chapter  IV of \cite{[LSU]},
the equation \eqref{7-eq-v} admits a unique solution $
v \in W_p^{2,1}( Q_{t_5 } ) $, and by virtue of
an extension of the same theorem
(see the last remark before paragraph 10, chapter IV of
\cite{[LSU]}) we obtain a unique solution
$ T \in W^{2,1}_q(Q_{ t_5})$ of the equation
\eqref{7-eq-T}.

\noindent Taking into account Lemma 6.3, we obtain (see \cite{[So76]} and
also Lemma 3.4 of the chapter II of \cite{[LSU]}) for $ 0<t \leq
t_6 $
\begin{equation}\label{ex(7.11)}
V_{(p,v)}(t) \leq c ( \| v_0 \|_{ W^{2-{2\over p}}_p(\Omega)}^p +
\| F_v \|_{L^p(Q_t ) }^p )  ,
\end{equation}
\begin{equation}\label{ex(7.11)-T}
V_{(q,T)}(t) \leq c ( \| T_0 \|_{ W^{2-{2\over q}}_q(\Omega)}^q +
\| F_T \|_{L^q(Q_t ) }^q
),
\end{equation}
where $F_v $ and $F_T $ are respectively the second member of
\eqref{7-eq-v} and that of \eqref{7-eq-T}.

\noindent It is not difficult to see that
$$
\Big\| R_0 \nabla ( ( {\varrho \over \mu_a} + {\pi \over \mu_h} )
\overline T ) + \Big[ \int_0^\infty ( \sigma + \nu )dm + \varrho
+ \pi \Big] \nabla \Phi \Big\|_{L^p(\Omega)}^p \leq c (1 +
V_{(q,\overline T)}(t) \| \overline T \|_{W^2_q(\Omega)}^{p-q}) ,
$$
$$
\int_0^t \big\| \eta \sum_{i,j=1}^3 \Big( {\partial \overline{ v}_i \over
\partial x_j} + {\partial \overline {v}_j \over \partial x_i} - {2 \over
3} \delta_{ij} \nabla \cdot \overline v \Big) {\partial \overline{ v}_i \over
\partial x_j} + \zeta (\nabla \cdot \overline v)^2
\big\|_{L^q(\Omega)}^q dt' \leq  c \int_0^t \| \overline v
\|_{W^2_p(\Omega)}^{2q - p} dt' V_{(p, \overline v )} (t) .
$$
The other terms of $F_v$ and $F_T$ can be estimated in the usual way
(see also \eqref{H-gl}). Thus we deduce
\eqref{ex(7.8)}--\eqref{ex(7.10)} from
\eqref{ex(7.11)}--\eqref{ex(7.11)-T}. \ $\square$

\section{ Existence and uniqueness of the local solution}

To prove Theorem 3.1, we start with the following lemma.

\begin{lem}
There exist positive constants $\overline {R}_v$, $\overline{ R}_T$ and $ t_5
\in \, ]0,  t_4]$ such that, if $0< t \leq t_5$, $ \overline v \in
\Theta^{(v)}_t $, $ \overline T \in \Theta^{(T)}_t $ and if
$$
V_{(p,\overline v )}( t ) \leq  \overline{ R}_v , \quad V_{(q,\overline T
)}( t ) \leq  \overline{ R}_T ,
$$
then the solution $(v,T)$ of the equations
\eqref{7-eq-v}--\eqref{7-eq-T} with the conditions
\eqref{cond-v}--\eqref{cond-T} satisfies the inequalities
$$
V_{(p, v)}( t ) \leq  \overline{ R}_v , \quad V_{(q, T)}( t ) \leq  \overline
{R}_T .
$$
\end{lem}

{\bf Proof.} \ The lemma follows from \eqref{ex(7.8)}--\eqref{ex(7.10)} by simple calculations.
\ $\square$

\medskip

We define
\begin{equation}\label{ex(8.1)} B_{t } =\{\, (v,T) \in
\Theta^{(v)}_{ t } \times \Theta^{(T)}_{ t } \, | \, V_{( p, v)}(t
) \leq  \overline{ R}_v , \ V_{( q, T)}(t ) \leq \overline{ R}_T \,  \} .
\end{equation}

For $0< t \leq t_5$ we define the map
$G_t:B_t\rightarrow \Theta^{(v)}_{ t } \times \Theta^{(T)}_{ t } $
such that, for $(\overline v, \overline T) \in B_{t }$,
$(v,T)=G_t(\overline v, \overline T)$ is the solution of the equation
\eqref{7-eq-v}--\eqref{7-eq-T} with the
conditions \eqref{cond-v}--\eqref{cond-T}.
By virtue of Lemma 8.1 we have
\begin{equation}\label{ex(8.2)}
G_t(B_{t })\subseteq B_{t }, \qquad 0< t \leq t_5 .
\end{equation}

\medskip
\medskip

{\bf Proof of Theorem 3.1.}  \ For $ 0< t \leq t_5 $, we define
\begin{equation}\label{ex(8.3)}
Y_t = \big[ L^2(0, t ; H^1_0 (\Omega ) ) \cap L^\infty (0, t
;L^2(\Omega)) \big] \times \big[ L^2(0, t ; H^1_0 (\Omega ) ) \cap
L^\infty (0, t ;L^2(\Omega)) \big] .
\end{equation}
We remark that the set $B_{t } $ defined in \eqref{ex(8.1)} is a
closed convex set in the space $Y_t $. Therefore, to prove the theorem,
it is sufficient to check that there exists $\overline t \in \; ]0,
t_5]$ such that the operator $G_{\bar t} $ is a contraction in the
natural topology of $Y_{\overline t} $.

Let be $ (\overline v_1, \overline T_1 ) $, $ (\overline v_2,
\overline T_2 ) \in B_{ t}$, $0< t \leq t_5$. First we consider
the solutions $( \varrho_1, \pi_1 , \sigma_1)$ and $( \varrho_2,
\pi_2, \sigma_2) $ for the equation system
\eqref{eq-varrho-lin}--\eqref{eq-sigma-lin} with the conditions
\eqref{cond-varrho}--\eqref{cond-sigma} and with the
substitutions $\overline v = \overline v_1$, $\overline T =
\overline T_1 $ and $\overline v = \overline v_2 $, $\overline T =
\overline T_2 $. We put
$$
E^{[\varrho]} = \varrho_1 - \varrho_2 , \qquad E^{[\pi]} = \pi_1 -
\pi_2 , \qquad E^{[\sigma]} = \sigma_1 - \sigma_2 ,
$$
$$
\overline D^{[ v]} = \overline v_1 - \overline v_2 , \qquad
\overline D^{[ T]} = \overline T_1 - \overline T_2 .
$$
We define
$$
{\overline u}_i ( t,m,x) = {\overline v}_i (t, x) -
\frac{1}{\alpha_l ( m )} \nabla \Phi (x) , \quad *\overline U_{4l,i} =
( mh_{gl} ( \overline T_i ,\pi_{i,\vartheta }
; m ), \overline u_{i,1}, \overline u_{i,2}, \overline u_{i,3})^T
,
$$
$$
\mathcal{E}^{[i]}=(\mathcal{E}^{[i]}_1,\mathcal{E}^{[i]}_2,\mathcal{E}^{[i]}_3), \quad i=1,2.
$$

From the difference between the equations
\eqref{eq-varrho-lin}-\eqref{eq-sigma-lin} for $ \varrho_1 ,\pi_1
,\sigma_1 $ and those for $\varrho_2 ,\pi_2 ,\sigma_2 $ it follows that
\begin{equation}\label{ex(8.6-rho)}
\partial_t E^{[\varrho]} +
\overline v_1 \cdot \nabla E^{[\varrho]} + \overline D^{[ v]}
\cdot \nabla \varrho_2 + E^{[\varrho]} \nabla \cdot \overline v_1
+ \varrho_2 \nabla \cdot \overline D^{[ v]} = 0,
\end{equation}
\begin{equation}\label{ex(8.6-pi)}
\partial_t E^{[\pi]} +
\overline v_1 \cdot \nabla E^{[\pi]} + \overline D^{[ v]} \cdot
\nabla \pi_2 + E^{[\pi]} \nabla \cdot \overline v_1 + \pi_2 \nabla
\cdot \overline D^{[ v]} =
\end{equation}
$$
= H_{gl} ( \overline T_2 ,\pi_{2} ,\sigma_2 ) - H_{gl} ( \overline
T_1 ,\pi_{1 } ,\sigma_1 ),
$$
\begin{equation}\label{ex(8.6-sigma)}
\partial_t E^{[ \sigma ]} +
(\overline U_{4l,1} - \overline U_{4l,2} ) \cdot \nabla_{( m,x)}
E^{[ \sigma ]} + \overline D^{[ U_{4l}]} \cdot \nabla_{( m,x)}
\sigma_2  +
\end{equation}
$$
+E^{[ \sigma ]} \nabla_{(m,x)} \cdot \overline U_{4l,1 } +
\sigma_2 \nabla_{(m,x)} \cdot \overline D^{[ U_{4l} ]}  =
$$
$$
= \big[ h_{gl}( \overline T_1 ,\pi_{1} ;m ) + B_1 (\sigma_1 ;m) -
g_1 ( m ) [ \pi_1  - \overline \pi_{vs(l)} (\overline T_1 )]^-
\big] E^{[ \sigma ]}  +
$$
$$
+ \big\{ h_{gl} (\overline T_1 ,\pi_{1} ;m) - h_{gl} (\overline
T_2 ,\pi_{2} ;m) + B_1 ( \sigma_1 ;m) - B_1 ( \sigma_2 ;m) +
$$
$$
- g_1 ( m ) \big( [ \pi_1  - \overline \pi_{vs(l)} ( \overline T_1
)]^- - [ \pi_2 - \overline \pi_{vs(l)} (\overline T_2 )]^- \big)
\big\}\sigma_2 +
$$
$$
+ g_0( m )\big( [\pi_1 - \overline \pi_{vs(l)} ( \overline T_1
)]^+ [ N^* - N ( \sigma_1  ) ]^+  - [ \pi_2 - \overline
\pi_{vs(l)} ( \overline T_2 ) ]^+ [N^* - N ( \sigma_2  ) ]^+ \big)
+
$$
$$
+ B_2 ( \sigma_1 ;m) - B_2 ( \sigma_2 ;m).
$$

We remember the relations
$$
\int_{\Omega}{(\overline v_1 \cdot \nabla E^{[\varrho]} )
E^{[\varrho]} }dx = -{1 \over 2} \int_{\Omega} (\nabla \cdot
\overline v_1 ) (E^{[\varrho]} )^{2}dx ,
$$
$$
\int_{\Omega}{(\overline v_1 \cdot \nabla E^{[\pi]} ) E^{[\pi]}
}dx = -{1 \over 2} \int_{\Omega} (\nabla \cdot \overline v_1 )
(E^{[\pi]} )^{2}dx ,
$$
$$
\int_{\Omega_{\overline M_1}}{(\overline U_{4l,1 } \cdot \nabla E^{[\sigma]} )
E^{[\sigma]} }dx = -{1 \over 2} \int_{\Omega_{\overline M_1}} (\nabla \cdot \overline
U_{4l,1 } ) (E^{[\sigma]} )^{2}dx .
$$
So, multiplying \eqref{ex(8.6-rho)}--\eqref{ex(8.6-sigma)}
by $ E^{[\varrho]} $, $ E^{[\pi]} $, $ E^{[\sigma]} $, integrating
the first two on $\Omega$ and the others on  $\Omega_{\overline
M_1}$ and taking into account estimates already known,  we obtain with
usual calculations
\begin{equation}\label{ex(8.4-rho)}
{d \over dt} \|E^{[\varrho]} \|^2_{L^2 ( \Omega )} \leq c ( 1 + \|
\overline v_1 \|_{W_p^2 ( \Omega )} ) \| E^{[ \varrho]} \|_{L^2 (
\Omega )}^2  + c \|\overline D^{[ v ]} \|_{H^1 ( \Omega)}^2 ,
\end{equation}
\begin{equation}\label{ex(8.4-pi)}
{d \over dt} \|E^{[\pi]}\|^2_{L^2 ( \Omega )} \leq c (1 + \|
\overline v_1 \|_{W_p^2 ( \Omega)} + \| \overline T_1 \|_{W_q^2 (
\Omega )}^2  ) \times
\end{equation}
$$
\times ( \| E^{[ \pi ]} \|_{L^2 ( \Omega )}^2  + \|
E^{[\sigma]}\|_{L^2 ( \Omega_{\overline M_1} )}^2) + c ( \| \overline D^{[ v ]}
\|_{H^1 ( \Omega )}^2 + \|\overline D^{[T]} \|_{L^2 ( \Omega )}^2
) ,
$$
\begin{equation}\label{ex(8.4-sigma)}
{d \over dt} \|E^{[\sigma]} \|^2_{L^2 ( \Omega_{\overline M_1})} \leq c (1 + \|
\overline v_1 \|_{W_p^2 ( \Omega)} + \| \overline T_1 \|_{W_q^2 (
\Omega )}^2 + \| \overline T_2 \|_{W_q^2 ( \Omega )}^2  ) \times
\end{equation}
$$
\times ( \| E^{[ \pi ]} \|_{L^2 ( \Omega )}^2  + \|
E^{[\sigma]}\|_{L^2 ( \Omega_{\overline M_1} )}^2) + c ( \| \overline D^{[ v ]}
\|_{H^1 ( \Omega )}^2 + \|\overline D^{[T]} \|_{L^2 ( \Omega )}^2
) ,
$$
From \eqref{ex(8.4-rho)}--\eqref{ex(8.4-sigma)} and the initial
conditions
$$
E^{[ \varrho ]} ( 0,\cdot  ) = 0, \qquad E^{[ \pi ]} (0,\cdot ) =
0, \qquad E^{[ \sigma ]} (0, \cdot , \cdot ) = 0,
$$
it follows that
\begin{equation}\label{ex(8.7-densi)}
\|E^{[\varrho]} (t)\|^2_{L^2 ( \Omega )} + \| E^{[ \pi ]} (t)
\|_{L^2 ( \Omega )}^2  + \| E^{[\sigma]}(t) \|_{L^2 ( \Omega_{\overline M_1} )}^2
 \leq
\end{equation}
$$
\leq c e^{c [ t +  t^{\frac{{p - 1}}{p}} +  t^{\frac{{q - 2}}{q}}
] } \int_0^t ( \| \overline D^{[ v ]} (t')\|_{H^1 ( \Omega )}^2 +
\|\overline D^{[T]} (t')\|_{L^2 ( \Omega )}^2 ) dt' .
$$

We consider the difference between the equations
\eqref{7-eq-v}--\eqref{7-eq-T} for $ ( v_1 ,T_1 ) = G_t( \overline
v_1 ,\overline T_1 ) $ and those for $ ( v_2,T_2 )= G_t (
\overline v_2 ,\overline T_2 ) $, so that if we put 
$D^{[ v ]} = v_1 - v_2$ and $ D^{[ T]} = T_1 - T_2 $ we obtain
\begin{equation}\label{ex(8.10-v)}
(\varrho_1 + \pi_1) \partial_t D^{[v]} - \eta \Delta D^{[v]} -
(\zeta + \frac{\eta }{3} ) \nabla ( \nabla \cdot D^{[v]})  = -
(E^{[\varrho]} +E^{[\pi]} ) \partial_t v_2 +
\end{equation}
$$
- (\varrho_1 + \pi_1) (\overline v_1 \cdot \nabla ) \overline
D^{[v]} - (E^{[\varrho]} + E^{[\pi]} ) (\overline v_1 \cdot \nabla
) \overline v_2 - (\varrho_2 + \pi_2) (\overline D^{[v]} \cdot
\nabla ) \overline v_2 +
$$
$$
- R_0 \nabla \big( \big( \frac{E^{[\varrho]} }{ \mu_a} +
\frac{E^{[\pi]} }{\mu_h }\big) \overline T_1 \big) - R_0\nabla
\big( \big( \frac{\varrho_2}{\mu_a } + \frac{\pi_2}{\mu_h} \big)
\overline D^{[T]} \big) -\Big[ \int_0^\infty  E^{[ \sigma]}
 dm + E^{[ \varrho ]}
+ E^{[ \pi ]} \Big]\nabla \Phi ,
$$
\begin{equation}\label{ex(8.10-T)}
(\varrho_1 + \pi_1) c_v \partial_t D^{[T]} - \kappa \Delta D^{[T]}
= - ( E^{[\varrho]} + E^{[\pi]}) c_v \partial_t T_2 +
\end{equation}
$$
- ( \varrho_1 + \pi_1 ) c_v \sum_{i = 1}^3 \overline v_{1,i}
\frac{\partial \overline D^{[T]} }{ \partial x_i } - ( \varrho_1 +
\pi_1 ) c_v \sum_{i = 1}^3 \overline D_{i}^{[v]} \frac{\partial
\overline T_2 }{\partial x_i } - ( E^{[\varrho]} + E^{[\pi]} ) c_v
\sum_{i = 1}^3 \overline v_{2,i} \frac{\partial \overline T_2
}{\partial x_i }  +
$$
$$
- R_0 ( \frac{\varrho_1}{\mu _a} + \frac{ \pi_1 }{\mu _h}
)\overline T_1 \nabla \cdot \overline D^{[v]}  - R_0 (
\frac{\varrho_1}{\mu _a} + \frac{ \pi_1 }{\mu _h} ) \overline
D^{[T]} \nabla \cdot \overline v_2 - R_0 (
\frac{E^{[\varrho]}}{\mu _a} + \frac{E^{[\pi]} }{\mu _h} )
\overline T_2 \nabla \cdot \overline v_2 +
$$
$$
+ \eta \sum_{i,j = 1}^3 \Big( \frac{ \partial \overline
D_{i}^{[v]}}{\partial x_j } + \frac{ \partial \overline
D_{j}^{[v]} }{ \partial x_i } - \frac{2}{3}\delta_{ij} \nabla
\cdot \overline D^{[v] } \Big) \frac{ \partial \overline v_{1,i}
}{\partial x_j } + \eta \sum_{i,j = 1}^3 \Big( \frac{ \partial
\overline v_{2,i} }{\partial x_j } + \frac{ \partial \overline
v_{2,j} }{\partial x_i } - \frac{2}{3}\delta _{ij} \nabla \cdot
\overline v_2 \Big) \frac{ \partial \overline D_{i}^{[v]}
}{\partial x_j } +
$$
$$
+ \zeta [ ( \nabla \cdot \overline v_1 )^2 - ( \nabla \cdot
\overline v_2 )^2 ] -\nabla \cdot
(\mathcal{E}^{[1]}-\mathcal{E}^{[2]}) + L_{gl} [ H_{gl} (
\overline T_1 , \pi_{1 } , \sigma_1 ) - H_{gl} ( \overline T_2 ,
\pi_{2} , \sigma_2 )] .
$$
If we remember \eqref{EQ4}, \eqref{def-E}
\eqref{cond0-P}, we have
\begin{equation}\label{def-div-E}
\nabla \cdot
(\mathcal{E}^{[1]}-\mathcal{E}^{[2]})=-\int_{0}^{\infty}
\big((a_{\lambda}^{(1)}+r_{\lambda}^{(1)})
\varrho_1+(a_{\lambda}^{(2)}+r_{\lambda}^{(2)})\pi_1+
\end{equation}
$$
+\int_0^{\infty}(a_{\lambda}^{(3)}(m)+r_{\lambda}^{(3)}(m))
\sigma_1(m)dm\big)\int_{S^2}(I^{(1)}_{\lambda}(x,q_1)-I^{(2)}_{\lambda}(x,q_1))dq_1d\lambda+
$$
$$
-\int_{0}^{\infty}
\big((a_{\lambda}^{(1)}+r_{\lambda}^{(1)}) E^{[\varrho]}+
(a_{\lambda}^{(2)}+r_{\lambda}^{(2)}) E^{[\pi]}+
$$
$$
+\int_0^{\infty}(a_{\lambda}^{(3)}(m)+
r_{\lambda}^{(3)}(m)) E^{[\sigma]}dm\big)\int_{S^2}
I^{(2)}_{\lambda}(x,q_1)dq_1d\lambda+
$$
$$
+\int_{0}^{\infty}\big(r_{\lambda}^{(1)}\varrho_1+r_{\lambda}^{(2)}\pi_1
+\int_{0}^{\infty}r_{\lambda}^{(3)}\sigma_1(m)dm\big)\int_{S^2}
(I^{(1)}_{\lambda}(x,q'_1)-I^{(2)}_{\lambda}(x,q'_1))dq'_1d\lambda
$$
$$
+\int_{0}^{\infty}\big(r_{\lambda}^{(1)}
 E^{[\varrho]}+r_{\lambda}^{(2)} E^{[\pi]}
+\int_{0}^{\infty}r_{\lambda}^{(3)}(m)
E^{[\sigma]}dm\big)\int_{S^2}
I^{(2)}_{\lambda}(x,q'_1)dq'_1d\lambda+
$$
$$
+4\pi\int_{0}^{\infty}
\big(a_{\lambda}^{(1)}\varrho_1+a_{\lambda}^{(2)}\pi_1
+\int_0^{\infty}a_{\lambda}^{(3)}(m)
\sigma_1(m)dm\big)(B(\lambda,\overline T_1)-B(\lambda,\overline
T_2))d\lambda+
$$
$$
+4\pi\int_{0}^{\infty}\big(a_{\lambda}^{(1)}
 E^{[\varrho]}+a_{\lambda}^{(2)} E^{[\pi]}
+\int_{0}^{\infty}a_{\lambda}^{(3)}(m)
E^{[\sigma]}dm\big)B(\lambda,\overline T_2)d\lambda.
$$

We multiply the equations \eqref{ex(8.10-v)} and
\eqref{ex(8.10-T)} respectively by $ \frac{D^{[ v ]} }{\varrho_1 +
\pi_1 }$ and $ \frac{D^{[T]} }{\varrho_1 + \pi_1 }$ and integrate
them on $\Omega $. We remember the inequalities
$$
\eta \Big| \int_\Omega \frac{ \sum_{i,j}^3 \frac{\partial
(\varrho_1 + \pi_1 )}{\partial x_j} \frac{\partial D^{[v]}_{i}
}{\partial x_j} D^{[v]}_{i} }{(\varrho_1 + \pi_1)^2} dx \Big| +
(\zeta + \frac{\eta}{3})\Big| \int_\Omega \frac{ \sum_{i,j}^3
\frac{\partial (\varrho_1 + \pi_1 )}{\partial x_i} \frac{\partial
D^{[v]}_{j} }{\partial x_j} D^{[v]}_{i} }{(\varrho_1 + \pi_1)^2}
dx \Big| \leq
$$
$$
\leq \overline c (\| \varrho_1\|_{W^1_p(\Omega)} + \|
\pi_1\|_{W^1_p(\Omega)} ) \|D^{[v]}
\|_{{H}^1(\Omega)}^{1+\frac{3}{p}} \| D^{[v]}
\|_{L^2(\Omega)}^{1-\frac{3}{p}} \leq
$$
$$
\leq \varepsilon \|D^{[v]} \|_{{H}^1(\Omega)}^2 + C_\varepsilon (
\| \varrho_1 \|_{W^1_p(\Omega)}^{\frac{2p}{p-3}}  + \| \pi_1
\|_{W^1_p(\Omega)}^{\frac{2p}{p-3}} ) \| D^{[v]}
\|^2_{L^2(\Omega)} ,
$$
where $C_\varepsilon$ \ is a constant determined by an arbitrary
constant $\varepsilon >0$. To estimate the term
$$ \Big|\int\limits_\Omega  {\frac{{D^{[T]}
}}{{\rho _1  + \pi _1 }}\nabla \cdot
(\mathcal{E}^{[1]}-\mathcal{E}^{[2]})dx} \Big| \le \frac{2}{{\inf
\rho _0 }}\left\| {\nabla  \cdot
(\mathcal{E}^{[1]}-\mathcal{E}^{[2]})} \right\|_{L^2 (\Omega )}
\left\| {D^{[T]} } \right\|_{L^2 (\Omega )},
$$
we deduce, from Lemma \ref{lemme-9-3}
\begin{equation}
\|\nabla\cdot
(\mathcal{E}^{[1]}-\mathcal{E}^{[2]})\|_{L^2(\Omega)}\leq
c(1+\|\varrho_1\|_{W^1_p(\Omega)}+
\|\pi_1\|_{W^1_p(\Omega)}+\|\sigma_1\|_{W^1_p(\Omega_{\overline M_1})})\times
\end{equation}
$$
\times( \|E^{[\varrho]} \|_{L^2 ( \Omega )} + \|
E^{[ \pi ]} \|_{L^2 ( \Omega )}  + \| E^{[\sigma]} \|_{L^2 ( \Omega_{\overline M_1}
)} + \|\overline D^{[T]} \|_{L^2 ( \Omega )} ).
$$

Thus using repeatedly the Sobolev,  H\"older and Cauchy-Schwartz
inequalities, we obtain
\begin{equation}\label{ex(8.11-v)}
{d \over dt}\| D^{[v]} \|^2_{L^2(\Omega)} + \overline c_0 \|D^{[v]}
\|^2_{{H}^1(\Omega)} \leq
\end{equation}
$$
\leq c ( 1 + \| \overline v_2 \|_{W_p^2 ( \Omega )^2 } ) ( \big\|
D^{[ v ]} \big\|_{L^2 ( \Omega  )}^2  + \big\| E^{[ \varrho ]}
\big\|_{L^2 ( \Omega )}^2  +
$$
$$
+ \big\| E^{[ \pi ]} \big\|_{L^2 ( \Omega )}^2 + \big\| E^{[
\sigma ]} \big\|_{L^2 ( \Omega_{\overline M_1} )}^2  + \big\| \overline D^{ [ v ]}
\big\|_{L^2 ( \Omega )}^2 + \big\| \overline D^{ [ T ]}
\big\|_{L^2 ( \Omega )}^2 ) ,
$$

\begin{equation}\label{ex(8.11-T)}
{d \over dt}\| D^{[T]} \|^2_{L^2(\Omega)} + \overline c_0 \|D^{[T]}
\|^2_{{H}^1(\Omega)} \leq
\end{equation}
$$
\leq c ( 1 + \| T _2 \|_{W_q^2 ( \Omega )}^2 + \| \overline v_2
\|_{W_p^2 ( \Omega )}^2 ) ( \big\| D^{[ T ]} \big\|_{L^2 ( \Omega
)}^2 + \big\| E^{[ \varrho ]} \big\|_{L^2 ( \Omega )}^2 +
$$
$$
+ \big\| E^{[ \pi ]} \big\|_{L^2 ( \Omega )}^2 + \big\| E^{ [
\sigma ]} \big\|_{L^2 ( \Omega_{\overline M_1} )}^2 + \big\| \overline D^{[ v ]}
\big\|_{L^2 ( \Omega )}^2 + \big\| \overline D^{[ T ]} \big\|_{L^2
( \Omega )}^2 ) .
$$

From  \eqref{ex(8.7-densi)}, \eqref{ex(8.11-v)}--\eqref{ex(8.11-T)}
we obtain
\begin{equation}\label{ex(8.20-bis)}
\| D^{[v]} (t)\|^2_{L^2(\Omega)} + \| D^{[T]}
(t)\|^2_{L^2(\Omega)} + \overline c_0 \int_0^t ( \|D^{[v]}
(t')\|^2_{{H}^1(\Omega)} + \|D^{[T]} (t')\|^2_{{H}^1(\Omega)} dt')
\leq
\end{equation}
$$
\leq c e^{c(t + t^{\frac{q-2}{q}}+t^{\frac{p-2}{p}})} (t +
t^{\frac{q-2}{q}}+t^{\frac{p-2}{p}}) \Big[ \| \overline
D^{[T]}\|_{L^\infty (0,t;L^2 ( \Omega))}^2 + \| \overline
D^{[v]}\|_{L^\infty (0,t;L^2 ( \Omega ))}^2 +
$$
$$
+ c e^{c [ t +  t^{\frac{{p - 1}}{p}} +  t^{\frac{{q - 2}}{q}} ] }
\int_0^t ( \| \overline D^{[ v ]} (t')\|_{H^1 ( \Omega )}^2 +
\|\overline D^{[T]} (t')\|_{L^2 ( \Omega )}^2 ) dt' \Big].
$$
The inequality \eqref{ex(8.20-bis)} allows us to find a $ \overline{t}\in [0,t_5]$ such that
$$
\| D^{[v]} \|^2_{L^\infty (0, \overline t; L^2(\Omega))} + \| D^{[T]}
\|^2_{L^\infty (0, \bar t; L^2(\Omega))} + \|D^{[v]} \|^2_{L^2 (0,
\overline t; {H}^1(\Omega))} + \|D^{[T]} \|^2_{L^2 (0, \overline t;
{H}^1(\Omega))} \leq
$$
$$
\leq \kappa \big( \| \overline D^{[v]} \|^2_{L^\infty (0, \overline t;
L^2(\Omega))} + \| \overline D^{[T]} \|^2_{L^\infty (0, \overline t;
L^2(\Omega))} + \|\overline D^{[v]} \|^2_{L^2 (0, \overline t;
{H}^1(\Omega))} + \|\overline D^{[T]} \|^2_{L^2 (0, \overline t;
{H}^1(\Omega))} \big)
$$
with $0 < \kappa < 1$. It means that the operator $G_{\overline t} :
B_{\overline t} \to B_{\overline t}$ is a contraction. This, also thanks to the lemma 4.1,  allows us to
conclude the proof of the existence and uniqueness of the solution
on the interval $[0,\overline t]$. \ $\square$

\end{document}